\documentclass[a4paper,11pt,reqno]{amsart}
\usepackage[centering,includeheadfoot,margin=2.5 cm]{geometry}
\usepackage{graphics,enumitem,epsfig,textcomp}
\usepackage{amsfonts,amsmath,amssymb,euscript,color,mathrsfs}
\usepackage[utf8]{inputenc}	

\usepackage[caption=false,subrefformat=parens,labelformat=parens]{subfig}

\usepackage{url}

\usepackage[symbol]{footmisc}

\numberwithin{equation}{section}
\newtheorem{theorem}{Theorem}

\newtheorem{lemma}{Lemma}
\newtheorem{proposition}{Proposition}
\newtheorem{remark}{Remark}

\def\d{\,\mathrm{d}}
\def\eps{\varepsilon}
\def\N{\mathbb{N}}
\def\R{\mathbb{R}}
\def\C{\hbox{\rlap{\kern.24em\raise.1ex\hbox
      {\vrule height1.3ex width.9pt}}C}}
\def\P{\hbox{\rlap{I}\kern.16em P}}
\def\Q{\hbox{\rlap{\kern.24em\raise.1ex\hbox
      {\vrule height1.3ex width.9pt}}Q}}
\def\M{\hbox{\rlap{I}\kern.16em\rlap{I}M}}
\def\Z{\hbox{\rlap{Z}\kern.20em Z}}
\def\({\begin{eqnarray}}
\def\){\end{eqnarray}}
\def\[{\begin{eqnarray*}}
\def\]{\end{eqnarray*}}
\def\part#1#2{\frac{\partial #1}{\partial #2}}

\def\grad{\nabla}
\def\Norm#1{\left\| #1 \right\|}
\def\pmb#1{\setbox0=\hbox{$#1$}
  \kern-.025em\copy0\kern-\wd0
  \kern-.05em\copy0\kern-\wd0
  \kern-.025em\raise.0433em\box0 }
\def\bar{\overline}

\def\tot#1#2{\frac{\d #1}{\d #2}} 
\def\laplace{\Delta}
\def\d{\,\mathrm{d}}
\def\N{\mathbb{N}}
\def\R{\mathbb{R}}

\def\epsilon{\varepsilon}

\def\E{\mathcal{E}}
\def\P{\mathbb{P}}
\def\Q{\mathbb{Q}}

\def\comment#1{\textcolor{blue}{\bf [#1]}}

\newcommand*\di{\mathop{}\!\mathrm{d}}
\newcommand\bl{\left(}
\newcommand\br{\right)}

\newcommand\Vset{V}
\newcommand\Eset{E}
\newcommand\adj{\text{adj}}
\newcommand\EsetT{E_T}
\newcommand\EsetS{E_S}
\newcommand\sgn{\text{sgn}}

\begin{document}


\centerline{{\Large\textbf{ODE and PDE based modeling of biological transportation networks}}}
\vskip 7mm

\centerline{
	{\large Jan Haskovec}\footnote{Mathematical and Computer Sciences and Engineering Division,
		King Abdullah University of Science and Technology,
		Thuwal 23955-6900, Kingdom of Saudi Arabia; 
		{\it jan.haskovec@kaust.edu.sa}}\qquad
	{\large Lisa Maria Kreusser}\footnote{Department of Applied Mathematics and Theoretical Physics (DAMTP), University of Cambridge, Wilberforce Road, Cambridge CB3 0WA, UK;
		{\it L.M.Kreusser@damtp.cam.ac.uk}}\qquad
	{\large Peter Markowich}\footnote{Mathematical and Computer Sciences and Engineering Division,
		King Abdullah University of Science and Technology,
		Thuwal 23955-6900, Kingdom of Saudi Arabia;
		{\it peter.markowich@kaust.edu.sa}}
	}
\vskip 10mm

\noindent{\bf Abstract.}
We study the global existence of solutions of a discrete (ODE based) model on a graph describing the formation of biological transportation networks, introduced by Hu and Cai. We propose an adaptation of this model so that a macroscopic (PDE based) system can be obtained as its formal continuum limit. We prove the global existence of weak solutions of the macroscopic PDE model. Finally, we present results of numerical simulations of the discrete model, illustrating the convergence to steady states, their non-uniqueness as well as their dependence on initial data and model parameters.
\vskip 7mm


\noindent{\bf AMSC:} 35B32, 35B36, 35K55, 35Q92, 70F10, 92C42




\tableofcontents


\section{Introduction}\label{sec:Intro}

Transportation networks are ubiquitous in living systems such as leaf venation in plants, mammalian circulatory systems that convey nutrients to the
body through blood circulation, or neural networks that transport electric charge.
Understanding the development, function and adaptation of biologic transportation networks has been a long standing interest of
the scientific community \cite{motivationbio1,motivationbio4,motivationbio2,motivationbio3}.
Mathematical modeling of transportation networks is traditionally based on discrete frameworks, in particular mathematical graph theory and discrete energy optimization,
where the energy consumption of the network is minimized under the constraint of constant total material cost.
However, networks and circulation systems in living organisms are typically subject to continuous adaptation,
responding to various internal and external stimuli.
For instance, for blood circulation systems it is well known that throughout the life of humans and animals, blood vessel systems are
continuously adapting their structures to meet the changing metabolic demand of the tissue. In particular,
it has been observed in experiments that blood vessels can sense the wall shear stress and adapt their
diameters according to it \cite{huadaptation}. Consequently, for biological applications it is necessary to employ the
dynamic class of models.

Motivated by this observation, Hu and Cai \cite{hu} introduced a new approach to dynamic modeling of transportation networks.
They propose a purely local dynamic adaptation model based on mechanical laws, consisting
of  a system of ordinary differential equations (ODE) on a graph, coupled to a linear system of equations (Kirchhoff law).
In particular, the model responds only to local information and fluctuations in flow distributions can be naturally incorporated.
Global existence of solutions of the coupled ODE-algebraic system is not trivial and, to our best knowledge,
has not been proved so far. The first goal of this paper is to close this gap.

In contrast to the discrete modeling approach, models based on systems of partial differential equations (PDE) can be used to describe
formation and adaptation of transportation networks based on macroscopic (continuum) physical laws.
Hu and Cai proposed a PDE-based continuum model \cite{hulecture} which was subsequently studied in the series of papers
\cite{bionetworksmodeling, bionetworksbookchapter, bionetworksanalysis, bionetworksnotes}.
The continuum model consists of a parabolic reaction-diffusion equation for the conductivity field, constrained by a Poisson equation
for the pressure field. However, no connection between the discrete (ODE-based) and continuum (PDE-based)
models for biological transportation networks has been established so far.

The second goal of this paper is to provide a formal continuum limit of an extension of the Hu and Cai model \cite{hu}
on regular equidistant grids; the rigorous limit passage will be studied in a consequent paper \cite{continuumlimit}.
The resulting continuum energy functional is of the form
\begin{align}  \label{Econt}
   \E[c] = \int_{\Omega} \grad p\cdot c\grad p + \frac{\nu}{\gamma} |c|^\gamma  \d x,
\end{align}
with the \emph{metabolic constant} $\nu>0$ and \emph{metabolic exponent} $\gamma>0$.
The energy functional is defined on the set of nonnegative diagonal tensor fields $c=c(x)$ on $\R^d$,
\( \label{cTensor}
   c = \begin{pmatrix} c^1 & & \\ & \ddots & \\ & & c^d \end{pmatrix}.
\)
The symbol $|c|^\gamma$ is defined as $\sum_{k=1}^d \left| c^k \right|^\gamma$.
The scalar pressure $p=p(x)$ of the fluid within the network (porous medium) is subject to
the Poisson equation
\begin{align} \label{Pcont}
   -\grad\cdot(c\grad p) = S,
\end{align}
equipped with no-flux boundary condition, and the datum $S=S(x)$ represents the intensity of sources and sinks.
The formal $L^2$-gradient flow (local dynamic adaptation model) of the energy \eqref{Econt}
constrained by \eqref{Pcont} is of the form
\( \label{GFcont}
   \partial_t c^k = \left( \partial_{x_k} p\right)^2 - \nu \bigl| c^k\bigr|^{\gamma-2} c^k, \qquad k=1,\ldots,d,
\)
subject to homogeneous Dirichlet boundary conditions, and coupled to \eqref{Pcont}. Clearly, the system suffers from two drawbacks: first, the possible strong degeneracy of the Poisson equation \eqref{Pcont},
and, second, the fact that \eqref{GFcont} is merely a family of ODEs, parametrized by the spatial variable.
Therefore, we shall consider a regularization/extension of \eqref{Pcont}--\eqref{GFcont}, where the Poisson equation is of the form
\begin{align} \label{PcontR}
   -\grad\cdot((r\mathbb{I} + c)\grad p) = S,
\end{align}
where $r=r(x) \geq r_0 > 0$ is a prescribed function that models the isotropic background permeability of the medium,
and $\mathbb{I}\in\R^{d\times d}$ is the unit matrix.
The second drawback is addressed by equipping the transient system \eqref{GFcont} with a linear diffusive term modeling random fluctuations in the medium,
\begin{align}\label{GFcontR}
   \partial_t c^k = D^2 \laplace c^k + \left( \partial_{x_k} p\right)^2 - \nu \bigl| c^k\bigr|^{\gamma-2} c^k, \qquad k=1,\ldots,d,
\end{align}
subject to homogeneous Dirichlet boundary conditions, where $D^2>0$ is the constant diffusivity.
Let us note that the model \eqref{PcontR}--\eqref{GFcontR} is a variant of  the tensor-based model proposed by D. Hu,
restricted to the set of diagonal tensors \cite{hucaitensor}.

The third goal of this paper is to prove the global existence of weak solutions of the PDE system \eqref{PcontR}--\eqref{GFcontR}.
The proof shall rely on the fact that it is a formal $L^2$-gradient flow of the regularized energy functional
\begin{align}  \label{EcontR}
   \E[c] = \int_{\Omega} D^2 |\grad c|^2 + \grad p\cdot (r\mathbb{I} + c)\grad p + \frac{\nu}{\gamma} |c|^\gamma  \d x,
\end{align}
where the symbol $|\grad c|^2$ is defined as $\sum_{k=1}^d \left| \grad c^k \right|^2$.

This paper is organized as follows.
In Section \ref{sec:micro} we describe the discrete model \cite{hu} introduced by Hu and Cai, establish its gradient flow structure
and prove the global existence of solutions of the corresponding ODE system coupled to the Kirchhoff law (linear system of equations).
In Section \ref{sec:formallimit} we motivate an adaptation of the Hu-Cai model so that a continuum model can be obtained as its formal macroscopic limit.
We then derive the PDE system \eqref{Pcont}--\eqref{GFcont} as the formal gradient flow of the continuum energy \eqref{Econt}
and prove the global existence of solutions for $\gamma>1$.
Finally, results of numerical simulations of the discrete Hu-Cai model are presented in Section \ref{sec:numerics},
illustrating the convergence to steady states, their non-uniqueness as well as their dependence on initial data and model parameters.

\section{The microscopic model}\label{sec:micro}


In this section we describe the microscopic model introduced by Hu and Cai \cite{hu} and reformulated in \cite{bionetworksbookchapter}. Let $G=(\Vset,\Eset)$ be an undirected connected graph, consisting of a finite set of vertices $\Vset$ and a finite set of edges $\Eset$ where the number of vertices is denoted by  $n=|\Vset|$. We assume that any pair of vertices is connected by at most one edge and a vertex is not connected to itself by an edge. We denote the edge  between vertices $i\in\Vset$ and $j\in\Vset$ by $(i,j)\in \Eset$. Since the graph is undirected we refer by $(i,j)$ and $(j,i)$ to the same edge. For each edge $(i,j)\in\Eset$ of the graph $G$ we consider its length and its conductivity, denoted by  $L_{ij}=L_{ji}>0$ and $C_{ij}=C_{ji}\geq 0$, respectively. 
In the sequel, we assume that the lengths $L_{ij}>0$ are given as a datum and fixed for all $(i,j)\in\Eset$.
The conductivities $C_{ij}$ are subject to the energy optimization and adaptation process. We assume that initially all edges in $\Eset$ have strictly positive conductivities.
In each vertex $i\in\Vset$ we have the pressure $P_i\in\R$.
The pressure drop between vertices $i\in\Vset$ and $j\in\Vset$ connected by an edge $(i,j)\in\Eset$ is given by 
\begin{align}\label{eq:pressuredrop}
(\Delta P)_{ij}:=P_j-P_i.
\end{align}
Note that the pressure drop is antisymmetric, i.e., by definition,  $(\Delta P)_{ij}=-(\Delta P)_{ji}$.
The oriented flux (flow rate) from vertex $i\in\Vset$ to $j\in\Vset$ is denoted by $Q_{ij}$; again, we have $Q_{ij}=-Q_{ji}$.
For biological networks, the Reynolds number of the flow is typically small and the flow is predominantly in the laminar (Poiseuille) regime.
Then the flow rate between vertices $i\in\Vset$ and $j\in\Vset$ along edge $(i,j)\in\Eset$ is proportional to the conductance $C_{ij}$ and the pressure drop $(\Delta P)_{ij}=P_j-P_i$,
\begin{align}\label{eq:flowrate}
Q_{ij} := C_{ij}\frac{P_j-P_i}{L_{ij}}\qquad\text{for all~}(i,j)\in \Eset.
\end{align}
The local mass conservation in each vertex is expressed in terms of the Kirchhoff law
\begin{align}\label{eq:kirchhoff}
-\sum_{j\in N(i)} C_{ij}\frac{ P_j-P_i}{L_{ij}}=S_i\qquad \text{for all~}i\in \Vset.
\end{align}
Here $N(i)$ denotes the set of vertices connected to $i\in\Vset$ through an edge,
and $S=(S_i)_{i\in\Vset}$ is the prescribed strength of the flow source ($S_i>0$) or sink ($S_i>0$) at vertex $i$.
Clearly, a necessary condition for the solvability of \eqref{eq:kirchhoff} is the global mass conservation
\begin{align}\label{eq:conservationmass}
\sum_{i\in\Vset}S_i=0,
\end{align}
which we assume in the sequel.
Given the vector of conductivities $C=(C_{ij})_{(i,j)\in\Eset}$, the Kirchhoff law \eqref{eq:kirchhoff} is a linear system of equations for the vector of pressures $P=(P_i)_{i\in\Vset}$.
With the global mass conservation \eqref{eq:conservationmass}, the linear system \eqref{eq:kirchhoff} is solvable if and only if the graph with edge weights $C=(C_{ij})_{(i,j)\in\Eset}$ is connected \cite{bionetworksbookchapter}, where only edges with positive conductivities $C_{ij}>0$ are taken into account
(i.e., edges with zero conductivities are discarded).
Note that the solution is unique up to an additive constant.

Hu and Cai \cite{hu} propose an energy cost functional consisting of a pumping power term and a metabolic cost term.
According to the Joule's law, the power (kinetic energy) needed to pump material through an edge $(i,j)\in\Eset$ is proportional to the pressure drop $(\Delta P)_{ij}=P_j-P_i$
and the flow rate $Q_{ij}$ along the edge, i.e.,
\begin{align*}
  (\Delta P)_{ij}Q_{ij}=\frac{Q_{ij}^2}{C_{ij}}L_{ij}. 
\end{align*}
The metabolic cost of maintaining the edge is assumed proportional to its length $L_{ij}$ and a power of its conductivity $C_{ij}^{\gamma}$, with an exponent $\gamma>0$ of the network.
For instance, in blood vessels the metabolic cost is proportional to the cross-section area of the vessel \cite{motivationblood}. Modeling the blood flow by Hagen-Poiseuille's law,
the conductivity is proportional to the square of the cross-section area, implying $\gamma=1/2$ for blood vessel systems.
For models of leaf venation the material cost is proportional to the number of small tubes, which is
proportional to $C_{ij}$, and the metabolic cost is due to the effective loss of the photosynthetic power
at the area of the venation cells, which is proportional to $C_{ij}^{1/2}$. Consequently, the effective value
of $\gamma$ typically used in models of leaf venation lies between $1/2$ and $1$, \cite{hu}.
Consequently, the energy cost functional is given by
\begin{align}\label{eq:energydisc}
\tilde{\E}[C] := \sum_{(i,j)\in\Eset}\bl \frac{Q_{ij}[C]^2}{C_{ij}}+\frac{\nu}{\gamma} C_{ij}^{\gamma}\br L_{ij},
\end{align}
where $Q_{ij}[C]$ is given by \eqref{eq:flowrate} with pressures calculated from the Kirchhoff's law \eqref{eq:kirchhoff},
and $\nu>0$ is the so-called metabolic coefficient.
Note that  every edge of the graph $G$ is counted exactly once in the above sum.

To compute the gradient flow of the energy \eqref{eq:energydisc} constrained by  Kirchhoff's law \eqref{eq:kirchhoff},
we need the following result about the derivative of the pumping term with respect to the conductivities:

\begin{lemma}\label{lem:derivative}
	Let $Q_{ij}[C] = C_{ij} \frac{P_j-P_i}{L_{ij}}$ for all $(i,j)\in E$ as in \eqref{eq:flowrate}, where $P$ is a solution
	of the linear system \eqref{eq:kirchhoff} with a given vector of conductivities $C$.
	Then, for any fixed $(k,l)\in E$ we have
	\(  \label{eq:derivativeq}
	\part{}{C_{kl}} \sum_{(i,j)\in E} \frac{Q_{ij}[C]^2}{C_{ij}}L_{ij} = - \frac{Q_{kl}[C]^2}{C_{kl}^2}L_{kl}.
	\)
\end{lemma}
\begin{proof}
Since
\begin{align*}
\part{}{C_{kl}} \sum_{(i,j)\in E} \frac{Q_{ij}[C]^2}{C_{ij}}L_{ij}=- \frac{Q_{kl}[C]^2}{C_{kl}^2}L_{kl}+2\sum_{(i,j)\in E} \frac{Q_{ij}[C]}{C_{ij}}\frac{\partial Q_{ij}[C]}{\partial C_{kl}}L_{ij},
\end{align*}
it is sufficient to show that
\begin{align*}
\sum_{(i,j)\in E} \frac{Q_{ij}[C]}{C_{ij}}\frac{\partial Q_{ij}[C]}{\partial C_{kl}}L_{ij}=0.
\end{align*}
Let $\mathbb{A}=(\mathbb{A}_{ij})$ denote the adjacency matrix of the graph $G=(\Vset,\Eset)$, i.e.\ its  coefficients are defined by
\begin{align}\label{adjM}
\mathbb{A}_{ij}=\begin{cases}
0 & \mbox{if }(i,j)\notin \Eset, \\
1 & \mbox{if }(i,j)\in \Eset.
\end{cases}
\end{align}
Note that $G$ is an undirected graph, implying $\mathbb{A}_{ij}=\mathbb{A}_{ji}$.
Due to the symmetry of $C_{ij}$ and $L_{ij}$ and antisymmetry of $Q_{ij}$ we have
\begin{align*}
2\sum_{(i,j)\in E} \frac{Q_{ij}[C]}{C_{ij}}\frac{\partial Q_{ij}[C]}{\partial C_{kl}}L_{ij}
&=\sum_{i=1}^n\sum_{j=1}^n \mathbb{A}_{ij}\bl \frac{P_j-P_i}{L_{ij}}\frac{\partial Q_{ij}[C]}{\partial C_{kl}}\br L_{ij}
\\&=\sum_{j=1}^n P_j\sum_{i=1}^n \mathbb{A}_{ij} \frac{\partial Q_{ij}[C]}{\partial C_{kl}} -\sum_{i=1}^n P_i\sum_{j=1}^n \mathbb{A}_{ij}  \frac{\partial Q_{ij}[C]}{\partial C_{kl}} 
\\&=-2\sum_{i=1}^n P_i\sum_{j=1}^n \mathbb{A}_{ij} \frac{\partial Q_{ij}[C]}{\partial C_{kl}} 
\\&=-2\sum_{i=1}^n P_i\frac{\partial}{\partial C_{kl}} \sum_{j\in N(i)} Q_{ij}.
\end{align*}
By the definition of the flow rate $Q_{ij}$ in \eqref{eq:flowrate} and Kirchhoff's law \eqref{eq:kirchhoff} we have
\begin{align*}
\sum_{j\in N(i)}Q_{ij }=\sum_{j\in N(i)} C_{ij} \frac{P_j-P_i}{L_{ij}}=S_i,
\end{align*}
and since the sources/sinks $S_i$ are fixed, we conclude
\begin{align*}
\sum_{(i,j)\in E} \frac{Q_{ij}[C]}{C_{ij}}\frac{\partial Q_{ij}[C]}{\partial C_{kl}}L_{ij}=0
\end{align*}
as required.
\end{proof}

Using the result in \eqref{eq:derivativeq} in the above Lemma, it is easy to see that for $(i,j)\in E$ the derivative of the energy \eqref{eq:energydisc} is given by
\begin{align*}
\frac{\partial}{\partial C_{ij}} \tilde{\E}[C] = - \left( \frac{Q_{ij}[C]^2}{C_{ij}^2} - \nu C_{ij}^{\gamma-1} \right) L_{ij}.
\end{align*}
Therefore, the  gradient flow of \eqref{eq:energydisc} constrained by  Kirchhoff's law \eqref{eq:kirchhoff} with respect
to the Euclidean distance 
is given by the ODE system
\begin{align}\label{eq:odesystem}
   \tot{C_{ij}}{t} = \left( \frac{Q_{ij}[C]^2}{C_{ij}^2} - \nu C_{ij}^{\gamma-1} \right) L_{ij},
\end{align}
coupled to  Kirchhoff's law \eqref{eq:kirchhoff} via the definition of the flow rate \eqref{eq:flowrate}.

The general formulation of a gradient flow of the functional $\tilde{\E}$ is of the form
\begin{align*}
\frac{\di z}{\di t}=-\mathcal{K}[z]\tilde{\E}'[z]\qquad\text{or, equivalently,}\qquad \mathcal{G}[z]\frac{\di z}{\di t}=-\tilde{\E}'[z],
\end{align*}
where $\tilde{\E}'(z)$ is the Fr\'echet derivative of the energy functional $\tilde{\E}\colon\mathcal{Z}\to\R$
on the subset $\mathcal{Z}$ of a linear space and $z\in \mathcal{Z}$.
We denote the space of tangent vectors at a point $z\in\mathcal{Z}$ by $\mathcal{T}_z\mathcal{Z}$
and the space of cotangent vectors, i.e.\ the set of all linear functionals on $\mathcal{T}_z\mathcal{Z}$, by $\mathcal{T}_z^*\mathcal{Z}$.
Then, the derivative $\tilde{\E}'[z]$ is a cotangent vector and $\mathcal{G}[z],\mathcal{K}[z]$ are duality maps,
mapping tangents to cotangents and vice versa, i.e.\ $\mathcal{G}[z]\colon \mathcal{T}_z\mathcal{Z}\to \mathcal{T}_z^*\mathcal{Z}$
and $\mathcal{K}[z]\colon \mathcal{T}_z^*\mathcal{Z}\to \mathcal{T}_z\mathcal{Z}$, with $\mathcal{K}=\mathcal{G}^{-1}$.
See, e.g., \cite{Peletier} for details.

Based on this general formulation, we consider the gradient flow with respect to a weighted Euclidean distance
and  introduce a duality map resulting in  the ODE system of the form
\begin{align}\label{eq:GF-alpha}
   \tot{C_{ij}}{t} = \left( \frac{Q_{ij}[C]^2}{C_{ij}} - \nu C_{ij}^{\gamma} \right) C_{ij}^{\alpha-1} L_{ij},
\end{align}
with a fixed exponent $\alpha\in\R$, constrained by the Kirchhoff law \eqref{eq:kirchhoff}.
For modeling reasons (see \cite{hu} and the references therein) we require that the speed of metabolic decay
is an increasing function of the conductivity. Therefore, we impose $\alpha > 1-\gamma$.
In particular, the choice $\alpha = 2-\gamma$
leads to the system studied by Hu and Cai in \cite{hu}.
Note that for $\alpha > 1-\gamma$ the solution of \eqref{eq:GF-alpha} is nonnegative for nonnegative initial data.
Moreover, we have the dissipation of the energy \eqref{eq:energydisc} along the solutions of \eqref{eq:energydisc}, \eqref{eq:kirchhoff}, since
\begin{align}\label{eq:EC-dissipation}
   \tot{}{t} \tilde{\E}[C] = \tilde{\E}'[C]\cdot \tot{C}{t}=  - \sum_{(i,j)\in \tilde{\E}} \left( \frac{Q_{ij}[C]^2}{C_{ij}} - \nu C_{ij}^{\gamma} \right)^2  C_{ij}^{\alpha-2} L_{ij}^2 \leq 0.
\end{align}

\subsection{Global existence of solutions} 
We shall prove the global existence of solutions for the ODE system \eqref{eq:GF-alpha}
coupled to the Kirchhoff law \eqref{eq:kirchhoff} through the definition of the flow rate \eqref{eq:flowrate}.
We assume that the initial datum for $C=C_{ij}$ is such that the underlying graph
is connected, where only edges with positive conductivity $C_{ij}>0$ are taken into account
(i.e., edges with zero initial conductivity are discarded and removed from the graph).
This implies that the Kirchhoff law \eqref{eq:kirchhoff} is solvable for $t=0$
(uniquely up to an additive constant) for the pressures.
Depending on the values of the exponents $\alpha\in\R$, $\gamma>0$, we distinguish two cases:
\begin{itemize}
	\item
	If $\gamma+\alpha \geq 2$, then we have for all $(i,j)\in E$,
	\(   \label{case1}
	    \tot{C_{ij}}{t} \geq -\nu L_{ij} C_{ij}^{\gamma+\alpha-1}.
	\)
	Then, since the exponent $\gamma+\alpha -1 \geq 1$, the solutions of \eqref{eq:GF-alpha} remain positive for all $t>0$
	(recall that the initial datum is strictly positive for all $(i,j)\in E$).
	Consequently, the underlying graph remains connected
	and the Kirchhoff law \eqref{eq:kirchhoff} remains solvable for all times.
	Moreover, the terms $C_{ij}^{\gamma}$ and $Q_{ij}^2/C_{ij}$ 
	remain globally bounded due to the energy dissipation \eqref{eq:EC-dissipation}.
	Thus, the solution of the system \eqref{eq:GF-alpha}, \eqref{eq:kirchhoff} exists globally in time.
		
	\item
	If $0<\gamma+\alpha -1 < 1$, the solution may exist only locally in time
	and some of the conductivities $C_{ij}$ may vanish in finite time.
	Then the edges with $C_{ij}=0$ are discarded and the connectivity
	of the graph may be lost, which would make the Kirchhoff law \eqref{eq:kirchhoff}
	unsatisfiable. However, as we prove below, under a certain assumption
	on the source/sink term $S_{ij}$, this does not happen, i.e.,
	the Kirchhoff law remains solvable even after the eventual removal of the edge(s)
	with vanishing conductivity.
	Thus, the solution $C=C(t)$ can be extended past this time simply by solving
	a reduced ODE system with initial datum equal to the 'terminal' state with the respective edge(s) removed.
	
	We start by proving a result stating that if we divide the set of vertices $\Vset$ into two disjoint parts $\Vset_1$, $\Vset_2$
such that the sources/sinks $S_{ij}$ induce a net flux $\Delta S\neq 0$ between them, 
then a connection (i.e., at least one edge with positive conductivity) between $\Vset_1$ and $\Vset_2$ will be maintained along the solutions
of \eqref{eq:GF-alpha}, \eqref{eq:kirchhoff}.

\end{itemize}

\begin{lemma}\label{lem:connect}
	Let $\gamma>0$ and $0 < \gamma+\alpha -1 < 1$.
	Let the set of vertices $\Vset$ be the disjoint union $\Vset_1 \cup \Vset_2$ such that
	\begin{align} \label{Sass}
	   \Delta S := \sum_{j\in\Vset_1} S_j = - \sum_{j\in\Vset_2} S_j \neq 0.
	\end{align}
	Let $\widetilde{\Eset}$ be the set of edges connecting $\Vset_1$ to $\Vset_2$,
	i.e.,
	\begin{align*}
	\widetilde{\Eset} = \{ (i,j)\in\Eset; \; i\in\Vset_1\,, j\in\Vset_2 \},
	\end{align*}
	and assume that $C_{ij}(t=0) \geq 0$ for all $(i,j)\in \widetilde{\Eset}$ with
	\begin{align} \label{lem_connect_ass}
	\sum_{(i,j)\in \widetilde{\Eset}} C_{ij}(t=0) > 0.
	\end{align}
	Then
	\begin{align}\label{connect_claim}
	\sum_{(i,j)\in \widetilde{\Eset}} C_{ij}(t) > 0\qquad\mbox{for all } t>0
	\end{align}
	along the solutions of \eqref{eq:GF-alpha}, \eqref{eq:kirchhoff}.
\end{lemma}

\begin{proof}
	For contradiction, assume that there exists a $T>0$ such that \eqref{connect_claim}
	holds for $t<T$ and
	\begin{align}\label{forContradiction}
	\lim_{t\to T-} \sum_{(i,j)\in \widetilde{\Eset}} C_{ij}(t) = 0.
	\end{align}
	For $t<T$ we have
	\begin{align}\label{eq:conn_ODE}
	\tot{}{t} \sum_{(i,j)\in \widetilde{\Eset}} C_{ij} = \sum_{(i,j)\in \widetilde{\Eset}} \left( Q_{ij}^2 C_{ij}^{\alpha-2} - \nu C_{ij}^{\gamma+\alpha-1} \right) L_{ij}.
	\end{align}
	Since $0 < \gamma+\alpha -1 < 1$, we have for $t<T$ and each $(k,l)\in \widetilde{\Eset}$ the inequality
	\begin{align*}
	C_{kl}^{\gamma+\alpha -1} \leq  \left( \sum_{(i,j)\in\widetilde{\Eset}} C_{ij} \right)^{\gamma+\alpha -1},
	\end{align*}
	where we used  that $C_{ij} \geq 0$ for $(i,j)\in\widetilde{\Eset}$ and $t<T$.
	Similarly, since $\alpha-2<0$, we have
	\begin{align*}
	C_{kl}^{\alpha-2} \geq  \left( \sum_{(i,j)\in\widetilde{\Eset}} C_{ij} \right)^{\alpha-2}.
	\end{align*}
	Inserting this into \eqref{eq:conn_ODE}, we obtain for $t<T$,
	\begin{align}\label{eq:conn_ODE2}
	\tot{}{t} \sum_{(i,j)\in \widetilde{\Eset}} C_{ij} \geq
	\left( \sum_{(i,j)\in\widetilde{\Eset}} C_{ij} \right)^{\alpha-2} \sum_{(i,j)\in \widetilde{\Eset}} Q_{ij}^2 L_{ij}
	- \nu \left( \sum_{(i,j)\in\widetilde{\Eset}} C_{ij} \right)^{\gamma+\alpha -1} \sum_{(i,j)\in\widetilde{\Eset}} L_{ij}.
	\end{align}
	
	Next, we shall estimate the term $ \sum_{(i,j)\in \widetilde{\Eset}} Q_{ij}^2 L_{ij}$ from below.
	Due to  Kirchhoff's law \eqref{eq:kirchhoffrescaled}, we have for $t<T$,
	\(   \label{DeltaS}
	\sum_{(i,j)\in \widetilde{\Eset}} Q_{ij} = \Delta S \neq 0.
	\)
	We claim that for each $t<T$ there exists an edge $(k,l)\in \widetilde{\Eset}$
	such that
	\[
	|Q_{kl}| \geq \frac{|\Delta S|}{ | \widetilde{\Eset} |}.
	\]
	If not, we would have
	\[
	\left| \sum_{(i,j)\in \widetilde{\Eset}} Q_{ij} \right| \leq \sum_{(i,j)\in \widetilde{\Eset}} |Q_{ij}| < |\Delta S|,
	\]
	a contradiction to \eqref{DeltaS}.
	Consequently, for each $t<T$ we estimate
	\[
 	   \sum_{(i,j)\in \widetilde{\Eset}} Q_{ij}^2 L_{ij}  \geq  \frac{|\Delta S|^2}{ | \widetilde{\Eset} |^2} \min_{(i,j)\in \widetilde{\Eset}} L_{ij}.
	\]
	Inserting this into \eqref{eq:conn_ODE2}, we obtain
	\(  \label{conn_ODE3}
	   \tot{}{t} u(t) \geq \kappa_1 u(t)^{\alpha-2} - \kappa_2 u(t)^{\gamma+\alpha -1},
	\)
	where we denoted
	\[
	   u(t) := \sum_{(i,j)\in \widetilde{\Eset}} C_{ij}(t),
	\]
	and the constants
	\[
	   \kappa_1 := \frac{|\Delta S|^2}{ | \widetilde{\Eset} |^2} \min_{(i,j)\in \widetilde{\Eset}} L_{ij} > 0,\qquad
	   \kappa_2 := \nu  \sum_{(i,j)\in\widetilde{\Eset}} L_{ij} > 0.
	\]
	Since according to the assumption \eqref{lem_connect_ass} we have $u(0)>0$, \eqref{conn_ODE3} implies that
	\[
	u(t) \geq \min \left\{ u(0), (\kappa_1/\kappa_2)^\frac{1}{\gamma+1}\right\} > 0
	\]
	for $t<T$, a contradiction to \eqref{forContradiction}.
\end{proof}

\begin{theorem}\label{thm:connect}
	Let $\gamma>0$ and $0 < \gamma+\alpha -1 < 1$.
	Assume that \eqref{Sass} holds for any disjoint sets $\Vset_1,\Vset_2\subset\Vset$ such that $\Vset=\Vset_1\cup\Vset_2$.
	Let the initial datum $C_{ij}(t=0)\geq 0$ be such that the graph induced by edges $(i,j)\in\Eset$ with $C_{ij}(t=0)>0$ is connected. 
	Then the  graph induced by the solutions $C_{ij}=C_{ij}(t)$ of \eqref{eq:GF-alpha}, \eqref{eq:kirchhoff}, where edges with vanishing conductivities are discarded,
 	remains connected for all times $t\geq 0$. In particular, solutions of \eqref{eq:GF-alpha}, \eqref{eq:kirchhoff} with removal of edges with vanishing conductivities
	exist globally in time.
\end{theorem}

\begin{proof}
Let us show that the graph remains connected for all times, i.e., for each $t>0$ there exist a path of edges with positive conductivity
connecting each pair of vertices. For contradiction, assume that at time $t_0>0$ no such path exists connecting a vertex $i\in\Vset$ to vertex $j\in\Vset$.
Then, collect all vertices connected by a path to $i\in\Vset$ in the set $\Vset_1$, and let $\Vset_2:=\Vset\setminus\Vset_1$ be its complement.
Since $i\in\Vset$ is not connected to $j\in\Vset$ at time $t_0$, also $\Vset_1$ is not connected to $\Vset_2$,
which is a contradiction to the statement of Lemma \ref{lem:connect}.

Consequently, the graph induced by the solutions $C_{ij}=C_{ij}(t)$ of \eqref{eq:GF-alpha}, \eqref{eq:kirchhoff} never becomes disconnected,
and thus, by the fundamental result of the graph theory \cite{graphtheory}, the Kirchhoff law \eqref{eq:kirchhoff} is solvable.
Moreover, since the terms $C_{ij}^{\gamma}$ remain globally bounded due to the energy dissipation \eqref{eq:EC-dissipation},
the solution does not blow up. It can only happen that some $C_{ij}$ vanish in finite time. In this case the corresponding edge(s)
are removed and the solution $C=C(t)$ is continued by solving 	a reduced ODE system.
In this way a global solution of the system \eqref{eq:GF-alpha}, \eqref{eq:kirchhoff} is constructed.
\end{proof}

\begin{remark}
The assumption of Theorem \ref{thm:connect} that \eqref{Sass} holds for any disjoint sets $\Vset_1,\Vset_2\subset\Vset$ such that $\Vset=\Vset_1\cup\Vset_2$
means that the graph cannot be partitioned into subgraphs with balanced sources/sinks (i.e., $\sum S_{ij}=0$ over the subgraph).
If the opposite is true, then the ODE system \eqref{eq:GF-alpha}, \eqref{eq:kirchhoff} can be solved separately for each of the subgraphs
(after eventual removal of edges connecting them).
\end{remark}

\section{The macroscopic model}\label{sec:formallimit}

The goal of this section is to derive the formal macroscopic limit of the discrete model \eqref{eq:kirchhoff}, \eqref{eq:energydisc} as the number of nodes and edges tends to infinity,
and to study the existence of weak solutions of the corresponding gradient flow.
The limit consists of an integral-type energy functional coupled to a Poisson equation.
We shall show that the derivation requires an appropriate rescaling of the Kirchhoff law \eqref{eq:kirchhoff} and of the energy functional \eqref{eq:energydisc}.
Moreover, we have to restrict ourselves to discrete graphs represented by regular grids, i.e.,  tessellation of the domain $\Omega\subset\R^d$, $d\in\N$,
by congruent identical parallelotopes.
This restriction is dictated by the requirement that the formal gradient flow of the rescaled energy functional, constrained by the rescaled Kirchhoff law,
is of the form \eqref{eq:odesystem}. 

\subsection{Rescaling of the Kirchhoff law}\label{sec:poissonrescale}
Let us denote the vertices left and right of vertex $i\in\Vset$ along the $k$-th spatial dimension by $(i-1)_k$ and, resp., $(i+1)_k$.
The Kirchhoff law \eqref{eq:kirchhoff} is than written as
\begin{align}\label{eq:kirchhofftomacro}
-\sum_{k=1}^d \bl C_{i,(i+1)_k}\frac{P_{(i+1)_k}-P_i}{L_{i,(i+1)_k}} - C_{(i-1)_k, i}\frac{P_i-P_{(i-1)_k}}{L_{(i-1)_k, i}}\br =S_i\qquad \text{for all~}i\in \Vset.
\end{align}
Our goal is to identify the Kirchhoff law with a finite difference discretization of the Poisson equation \eqref{Pcont},
\begin{align}\label{eq:poissonhelp}
-\nabla \cdot\left(  c\nabla p\right)=S,
\end{align}
where $S=S(x)$ is a formal limit of the sequence of discrete sources/sinks $S_i$.
Clearly, for this the edge lengths in the left-hand side of \eqref{eq:kirchhofftomacro} have to appear quadratically in the denominator
instead of linearly. Alternatively, we can say that the sources/sinks $S_i$ in the right-hand side of \eqref{eq:kirchhofftomacro} have to be rescaled appropriately,
reflecting the fact that the edges of the graph are inherently one-dimensional structures.
A straightforward calculation reveals that a finite difference discretization of \eqref{eq:poissonhelp},
where $c=c(x)$ is an appropriate limit of the sequence of discrete conductivities,
is obtained if and only if
\begin{align*}
\frac{2}{L_{(i-1)_k,i}+L_{i,(i+1)_k}}=\frac{1}{L_{(i-1)_k,i}}=\frac{1}{L_{i,(i+1)_k}}
\end{align*}
for all $i\in\Vset$ and for all directions $k=1,\dots,d$.
Therefore, grid points must be equidistant in each spatial dimension, and we denote $h_k>0$ the grid spacing in the $k$-th dimension.
The discrete graph is thus identified with a tessellation of $\Omega$ by identical parallelotopes.
For simplicity, we restrict ourselves to work with rectangular parallelotopes (bricks) in the sequel,
with edges parallel to the axes.
A generalization of the result for parallelotopes instead will be given in Remark \ref{rem:paralellogram}.
The rescaled Kirchhoff law is then written as
\begin{align}\label{eq:kirchhoffrescaled}
   -\sum_{k=1}^d \frac{1}{h_k} \bl C_{i,(i+1)_k}\frac{P_{(i+1)_k}-P_i}{h_k}-C_{(i-1)_k, i}\frac{P_i-P_{(i-1)_k}}{h_k}\br =S_i\qquad \text{for all~}i\in \Vset.
\end{align}

\subsection{Rescaling of the discrete energy functional}
In order to obtain an integral-type functional in the macroscopic limit of the sequence
of discrete energy functionals \eqref{eq:energydisc}, they need to be properly rescaled
depending on the spatial dimension $d\in\N$.
In particular, \eqref{eq:energydisc} has to be replaced by
\begin{align} \label{eq:energydiscweightshelp}
   \tilde{\E}[C] = \sum_{(i,j)\in\Eset}\bl \frac{Q_{ij}[C]^2}{C_{ij}}+\frac{\nu}{\gamma} C_{ij}^{\gamma}\br W_{ij}^d,
\end{align}
where $W_{ij}$ are some (abstract) weights that scale linearly with the grid spacing.
Before we introduce the formal macroscopic limit of the rescaled discrete functional \eqref{eq:energydiscweightshelp}
constrained by the rescaled Kirchhoff law \eqref{eq:kirchhoffrescaled}, let us make the following
observation about the gradient flow \eqref{eq:energydiscweightshelp}--\eqref{eq:kirchhoffrescaled}.

\begin{proposition}\label{prop:GFw-discrete}
Consider the setting introduced in Section \ref{sec:poissonrescale} with the discrete graph realized as a rectangular tessellation of $\Omega\in\R^d$.
Then the formal gradient flow (with respect to the Euclidean distance)
of the energy functional \eqref{eq:energydiscweightshelp} constrained by the rescaled Kirchhoff law \eqref{eq:kirchhoffrescaled}
is of the type \eqref{eq:odesystem}, i.e.,
 \begin{align}\label{eq:gradientflowrescaled}
    \tot{C_{ij}}{t} = \bl \frac{Q_{ij}[C]^2}{C_{ij}^2}-\nu  C_{ij}^{\gamma-1} \br W_{ij}^d,
 \end{align}
if and only if all the weights $W_{ij}$ are equal.
\end{proposition} 

\begin{proof} 
Denoting the adjacency matrix \eqref{adjM} of the tessellation by $\mathbb{A}=(\mathbb{A}_{ij})$, we have for any edge $(l,m)\in\Eset$,
\begin{align*}
   \frac{\partial \tilde{\E}[C]}{\partial C_{lm}}
      =-\frac{Q_{lm}[C]^2}{C_{lm}^2}W_{lm}+\nu  C_{lm}^{\gamma-1}W_{lm}
         +\frac{1}{2}\sum_{i=1}^n\sum_{j=1}^n \mathbb{A}_{ij}\bl \frac{2Q_{ij}[C]}{C_{ij}}\frac{\partial Q_{ij}[C]}{\partial C_{lm}}\br W_{ij}^d.
\end{align*}
The last term of the right-hand side is equal to
\begin{align}\label{eq:gradientflowhelp}
\begin{split}
&\sum_{i=1}^n\sum_{j=1}^n \mathbb{A}_{ij}\bl \frac{P_j-P_i}{L_{ij}}\frac{\partial Q_{ij}[C]}{\partial C_{lm}}\br W^d_{ij}\\&=-\sum_{j=1}^n P_j\sum_{i=1}^n \mathbb{A}_{ij} \frac{\partial Q_{ji}[C]}{\partial C_{lm}} \frac{W_{ij}^d}{L_{ij}}-\sum_{i=1}^n P_i\sum_{j=1}^n \mathbb{A}_{ij}  \frac{\partial Q_{ij}[C]}{\partial C_{lm}} \frac{W_{ij}^d}{L_{ij}}
\\&= -2\sum_{i=1}^n P_i\sum_{j=1}^n \mathbb{A}_{ij}  \frac{\partial Q_{ij}[C]}{\partial C_{lm}} \frac{W_{ij}^d}{L_{ij}}.
\end{split}
\end{align}
Now note that the rescaled Kirchhoff law \eqref{eq:kirchhoffrescaled} is in terms of $Q_{ij}$, $L_{ij}$ written as
\begin{align*}
   - \sum_{j\in\Vset} \mathbb{A}_{ij} \frac{Q_{ij}}{L_{ij}} =S_i\qquad \text{for all~}i\in \Vset.
\end{align*}
Therefore, if (and only if) all the weights $W_{ij}$ are equal to the same value $W>0$, we have
\begin{align*}
   \sum_{i=1}^n P_i\sum_{j=1}^n \mathbb{A}_{ij}  \frac{\partial Q_{ij}[C]}{\partial C_{lm}} \frac{W_{ij}^d}{L_{ij}} =
   \sum_{i=1}^n P_i W^d \frac{\partial}{\partial C_{lm}} \left( \sum_{j=1}^n \mathbb{A}_{ij}  \frac{Q_{ij}}{L_{ij}} \right) = 0,
\end{align*}
and we obtain \eqref{eq:gradientflowrescaled} as the gradient flow.
\end{proof}

Note that for the grid consisting of a rectangular tessellation, the natural choice of the weight $W_{ij}\equiv W$ is 
\begin{align}\label{eq:weightsrec}
W^d=\prod_{k=1}^d h_k, 
\end{align}
i.e., the area of the rectangles for $d=2$ and the volume of the bricks for $d=3$.

\subsection{Formal derivation of the macroscopic model}
In this Section we shall show that the rescaled Kirchhoff law represents a finite difference
discretization of the Poisson equation \eqref{Pcont}, and that the discrete energy functional
\eqref{eq:energydiscweightshelp} with \eqref{eq:weightsrec} is an approximation (Riemann sum) of the integral-type functional
\eqref{Econt}. We shall work in the setting introduced above, i.e., the discrete graph is realized as a rectangular tessellation
of the rectangular domain $\Omega\in\R^d$.

Let us consider $p=p(x)$ a solution of the Poisson equation \eqref{Pcont},
\begin{align*}
   -\nabla \cdot\left(c\nabla p\right) = S,
\end{align*}
subject to the no-flux boundary condition on $\partial\Omega$.
Here $c=c(x)$ is a given diagonal permeability tensor field
\begin{align*}
   c = \begin{pmatrix} c^1 & & \\ & \ddots & \\ & & c^d \end{pmatrix},
\end{align*}
with the scalar nonnegative functions $c^k\in C(\Omega)$, $k=1,\dots,d$.
The density of sources/sinks $S=S(x)$ is given as a datum and satisfies the global mass balance
\[
   \int_\Omega S(x) \d x = 0.
\]
As already mentioned in Section \ref{sec:Intro}, existence of solutions of \eqref{Pcont} is not guaranteed
due to the possible strong degeneracy of the permeability tensor. However, as we are interested in a formal
derivation only, we assume that $p=p(x)$ exists as a strong solution of \eqref{Pcont}, i.e., is at least $C^2_b$ on $\Omega$.
Moreover, we assume that the elements of $c=c(x)$ are at least $C^1_b$ on $\Omega$.
Since $c=c(x)$ is diagonal, the left-hand side of \eqref{Pcont} can be rewritten as
\[
   -\nabla \cdot\left(c\nabla p\right) = - \sum_{k=1}^d \partial_{x_k} (c^k \partial_{x_k} p).
\]
Let $X_i\in\Omega$ be the physical location of the vertex $i\in\Vset$.
Denoting the flux $q^k := c^k \partial_{x_k} p$, a finite difference approximation of the term $\partial_{x_k} q^k$
at $x=X_i$ reads
\(  \label{findiff1}
   \partial_{x_k} q^k(X_i) \approx \frac{q^k(X_{(i+1/2)_k}) - q^k(X_{(i-1/2)_k})}{h_k} + \mathcal{O}(h_k),
\)
where $X_{(i+1/2)_k}$ and, resp., $X_{(i-1/2)_k}$ denotes the midpoint
of the edge connecting $X_i$ to its adjacent vertex to the right and, resp., to the left
in the $k$-th spatial direction.
A finite difference approximation of $q^k$ at $X_{(i+1/2)_k}$ reads
\(  \label{findiff2}
   q^k(X_{(i+1/2)_k}) = c^k(X_{(i+1/2)_k}) \frac{p(X_{(i+1)_k}) - p(X_{(i-1)_k})}{h_k} + \mathcal{O}(h_k),
\)
where $X_{(i+1)_k}$, resp., $X_{(i-1)_k}$ denotes the adjacent vertex of $X_i$
to the right and, resp., to the left in the $k$-th spatial direction.
We discretize $q^k(X_{(i-1/2)_k})$ analogously.
Putting \eqref{findiff1} and \eqref{findiff2} together and denoting
\( \label{notation}
  \begin{aligned}
   C_{i,(i\pm 1)_k} := c^k(X_{(i\pm 1/2)_k}), \qquad S_i := S(X_i),\\
   P_i := p(X_i), \qquad P_{(i\pm 1)_k} := p(X_{(i\pm 1)_k}),
   \end{aligned}
\)
we conclude that the rescaled Kirchhoff law \eqref{eq:kirchhoffrescaled} is a first order
finite difference approximation of the Poisson equation \eqref{Pcont}.

With the choice \eqref{eq:weightsrec} for the weight $W$, we have for $k=1,\dots,d$ and $c^k\in C^1_b(\Omega)$,
\[
   \int_\Omega |c^k|^\gamma \d x = W \sum_{i\in\Vset} \left| c^k(X_{(i +1/2)_k}) \right|^\gamma + \mathcal{O}(h_k).
\]
Moreover, we have
\[
   \int_\Omega c^k (\partial_{x_k} p)^2 \d x = W \sum_{i\in\Vset} c^k(X_{(i + 1/2)_k}) \left( \frac{p(X_{(i+1)_k}) - p(X_i)}{h_k} \right)^2 + \mathcal{O}(h_k).
\]
Therefore, noting that for the rectangular grid the energy functional \eqref{eq:energydiscweightshelp} can be rewritten as
\begin{align*}
   \tilde{\E}[C]= \frac12 \sum_{k=1}^d \sum_{i\in\Vset} \sum_{j\in N(i;k)} \bl \frac{Q_{ij}[C]^2}{C_{ij}}+\frac{\nu}{\gamma} C_{ij}^{\gamma}\br h_k^d,
\end{align*}
we have with the notation \eqref{notation}
\[
   \tilde{\E}[C] = {\E}[c] + \mathcal{O}(h),
\]
with the continuum energy defined by \eqref{Econt}, i.e.,
\[  \label{Econt2}
   \E[c] = \int_{\Omega} \grad p\cdot c\grad p + \frac{\nu}{\gamma} |c|^\gamma  \d x,
\]
where we recall that the symbol $|c|^\gamma$ is defined as $\sum_{k=1}^d \left| c^k \right|^\gamma$.

We now calculate the formal $L^2$-gradient flow of the energy \eqref{Econt} constrained by the Poisson equation \eqref{Pcont}.

\begin{lemma}\label{lem:gfcont}
The formal $L^2$-gradient flow of the continuum energy functional \eqref{Econt} constrained by the Poisson equation \eqref{Pcont} is given by
\eqref{GFcont}, i.e., 
\[
   \partial_t c^k = \left( \partial_{x_k} p\right)^2 - \nu \bigl| c^k \bigr|^{\gamma-2} c^k.
\]
\end{lemma}

\begin{proof}
Let us calculate the first variation of $\E$ in the direction $\phi$ where $\phi$ denotes a diagonal matrix with entries $\phi^1,\ldots,\phi^d$.
Using the expansion 
\begin{align}\label{eq:expansioncont}
   p[c+\epsilon\phi]=p_0+\epsilon p_1 + \mathcal{O}(\eps^2),
\end{align}
we have
\begin{align}\label{eq:firstvariation}
   \left.\tot{}{\epsilon}\E[c+\epsilon\phi]\right|_{\epsilon=0} =
      \sum_{k=1}^d\int_{\Omega} \left( \partial_{x_k} p_0\right)^2 \phi^k + 2 c^k (\partial_{x_k} p_0)(\partial_{x_k} p_1) + \nu \bigl| c^k\bigr|^{\gamma-2} c^k \phi^k \di x.
\end{align}
Multiplication the Poisson equation \eqref{Pcont} with permeability tensor $c+\epsilon \phi$ by $p_0$ and integration by parts gives
\begin{align*}
   \sum_{k=1}^d\int_{\Omega} \bl c^k+\epsilon\phi^k\br\left( \partial_{x_k} p_0\right)^2 + \epsilon c^k (\partial_{x_k} p_0) (\partial_{x_k} p_1) \di x
   = \int_{\Omega} Sp_0\di x + \mathcal{O}(\epsilon^2).
\end{align*}
Subtracting the identity
\[
   \sum_{k=1}^d \int_{\Omega} c^k \left( \partial_{x_k} p_0\right)^2 \di x= \int_{\Omega} S p_0\di x,
\]
we obtain
\begin{align*}
   \sum_{k=1}^d\int_{\Omega} \left( \partial_{x_k} p_0\right)^2 \phi^k + c^k (\partial_{x_k} p_0)(\partial_{x_k} p_1) \di x = 0.
\end{align*}
Plugging this into \eqref{eq:firstvariation} gives
\begin{align*}
   \left.\tot{}{\epsilon} \E[c+\epsilon\phi]\right|_{\epsilon=0} =
   \sum_{k=1}^d\int_{\Omega} \left[ - \left( \partial_{x_k} p_0\right)^2 + \nu \bigl| c^k\bigr|^{\gamma-2} c^k \right] \phi^k\di x.
\end{align*}
\end{proof}

\begin{remark}\label{rem:paralellogram}
We can easily generalize to the situation when the grid is 
realized by congruent identical parallelotopes with edges 
in linearly independent directions $\theta_1, \dots, \theta_d\in\R^d$.
Then the coordinate transform $e_k \mapsto \theta_k$ in \eqref{Econt}--\eqref{Pcont}, where $e_k$ is the
$k$-th vector of the generic basis of $\R^d$, leads to the transformed continuum energy functional
\begin{align}  \label{EcontT}
   E[c] = \int_{\Omega} \grad p\cdot \mathbb{P}[c]\grad p + \frac{\nu}{\gamma} |c|^\gamma  \d x,
\end{align}
coupled to the Poisson equation
\[
   -\grad\cdot \left(\mathbb{P}[c]\grad p \right) = S
\]
with the permeability tensor
\[
   \mathbb{P}[c] := \sum_{k=1}^d c^k\theta_k\otimes\theta_k.
\]
The corresponding formal $L^2$-gradient flow is of the form
\[
   \partial_t c^k = \left( \theta_k\cdot\grad p\right)^2 - \nu \bigl| c^k\bigr|^{\gamma-2} c^k.
\]
\end{remark}

\subsection{Global existence of solutions of a modified macroscopic model}
As noted in Section \ref{sec:Intro}, the model \eqref{Pcont}--\eqref{GFcont} suffers from two drawbacks:
first, the Poisson equation \eqref{Pcont} is possibly strongly degenerate since in general
the eigenvalues (i.e., diagonal elements) of the permeability tensor $c=c(x)$ may vanish.
To overcome this problem, we introduce a regularization of \eqref{Pcont} of the form
\begin{align} \label{PcontR2}
   -\grad\cdot(\P[c]\grad p) = S,
\end{align}
with the permeability tensor
\( \label{P}
   \P[c] := r\mathbb{I} + c,
\)
where $r=r(x) \geq r_0 > 0$ is a prescribed function that models the isotropic background permeability of the medium,
and $\mathbb{I}\in\R^{d\times d}$ is the unit matrix.
Clearly, \eqref{PcontR2} is uniformly elliptic as long as the eigenvalues of $c=c(x)$ are nonnegative.

The second drawback is due to the fact that \eqref{GFcont} is merely a family of ODEs, parametrized by the spatial variable $x\in\Omega$.
We cure this problem by introducing a linear diffusive term modeling random fluctuations in the medium. We thus obtain
\begin{align}\label{GFcontR2}
   \partial_t c^k = D^2 \laplace c^k + \left( \partial_{x_k} p\right)^2 - \nu \bigl| c^k\bigr|^{\gamma-2} c^k, \qquad k=1,\ldots,d,
\end{align}
subject to homogeneous Dirichlet boundary data, where $D^2>0$ is the constant diffusivity.
By a simple modification of the proof of Lemma \ref{lem:gfcont} we conclude that the system  \eqref{PcontR2}--\eqref{GFcontR2}
represents the formal $L^2$-gradient flow of the energy functional
\begin{align}  \label{EcontR2}
   \E[c] = \int_{\Omega} \frac{D^2}{2} |\grad c|^2 + \grad p\cdot \P[c]\grad p + \frac{\nu}{\gamma} |c|^\gamma  \d x,
\end{align}
with $\P[c]$ given by \eqref{P}, the symbol $|c|^\gamma$ is defined as $\sum_{k=1}^d \left| c^k \right|^\gamma$
and the symbol $|\grad c|^2$ is defined as $\sum_{k=1}^d \left| \grad c^k \right|^2$.
The gradient flow property is fundamental for proving the global existence of weak solutions of the PDE system \eqref{PcontR2}--\eqref{GFcontR2}.
We consider the PDE system on a bounded domain $\Omega\subset \R^d$ with smooth boundary $\partial \Omega$,
subject to homogeneous Dirichlet boundary conditions for $c$ and no-flux boundary conditions for  $p$, 
\begin{align}\label{eq:contbc}
   c(t,x)=0,\quad \frac{\partial p}{\partial n}(t,x)=0\quad \text{for }x\in\partial\Omega,\; t\geq 0,
\end{align}
where $n$ denotes the exterior normal vector to the boundary $\partial \Omega$. 
Moreover, we prescribe the initial datum for $c$,
\begin{align}\label{eq:continitialdata}
   c(t=0,x)=c^I(x)\quad \text{for }x\in\Omega,
\end{align}
where $c^I = c^I(x)$ is a diagonal tensor field in $\R^{d\times d}$ with nonnegative diagonal elements.

\begin{theorem}\label{th:existencecont}
	Let $S\in L^2(\Omega)$, $\gamma>1$ and $c^I\in H_0^1(\Omega)^{d\times d}\cap L^{\gamma}(\Omega)^{d\times d}$.
	Then the system \eqref{PcontR2}--\eqref{GFcontR2} subject to the data \eqref{eq:contbc}--\eqref{eq:continitialdata}
	admits a global weak solution $(c,p)$ such that
	\begin{align} \label{weaksol}
	\begin{split}
	   c\in L^{\infty}(0,\infty;H_0^1(\Omega))\cap L^{\infty}(0,\infty;L^{\gamma}(\Omega)),\quad \partial_t c\in L^2((0,\infty)\times\Omega),\\
	   \nabla p\in L^{\infty}(0,\infty;L^2(\Omega)),\quad c\nabla p\in L^{\infty}(0,\infty;L^2(\Omega)).
	\end{split}
	\end{align} 
	This solution satisfies the energy dissipation inequality
	\begin{align}\label{eq:energydissinequ}
	   \E[c(t)]+\sum_{k=1}^d \int_0^t \int_{\Omega} \bl \partial_t c^k(s,x)\br^2 \di x\di s\leq \E[c^I]\quad \text{for all~}t\geq 0,
	\end{align}
	with $\E[c]$ given by \eqref{EcontR2}.
\end{theorem}

For the proof of the above Theorem we adopt a strategy similar to \cite{bionetworksanalysis, bionetworksnotes}:
For $\eps>0$ we introduce the regularized Poisson equation
\begin{align}\label{eq:poissonregeta}
   -\grad\cdot(\P^\epsilon[c]\grad p) = S
\end{align}
with the permeability tensor
\( \label{Peps}
   \P^\eps[c]:= r\mathbb{I} + c\ast\eta_\eps,
\)
subject to no-flux boundary data for $p$.
Here, $\eta_\eps$ is a nonnegative, radially symmetric mollifier and the convolution $c\ast\eta_\eps$ is carried out elementwise,
\[
   c^k\ast\eta_\eps(x) := \int_{\R^d} c^k(y) \eta_\eps(x-y) \d y.
\]
Moreover, we regularize \eqref{GFcontR2} as follows,
\begin{align}\label{eq:gradientflowcontreg}
   \frac{\partial c^k}{\partial t}=D^2\Delta c^k+ \left( \partial_{x_k} p\right)^2\ast \eta_\eps - \nu \bigl| c^k \bigr|^{\gamma-2} c^k,\qquad k=1,\ldots,d.
\end{align}
By a slight adaptation of the proof of Lemma \ref{lem:gfcont} it is easily shown that \eqref{eq:poissonregeta}--\eqref{eq:gradientflowcontreg}
is the formal $L^2$-gradient flow of the energy
\begin{align}\label{eq:energycontregeta}
   \E^\eps[c] := \int_{\Omega} \frac{D^2}{2} \left|\grad c\right|^2 + \grad p\cdot \P^\epsilon[c]\grad p  + \frac{\nu}{\gamma} |c|^\gamma  \d x,
\end{align}
where we used the notation
\[
   |\grad c|^2 := \sum_{k=1}^d \left| \grad c^k \right|^2,\qquad
   |c|^\gamma := \sum_{k=1}^d \left| c^k \right|^\gamma.
\]

For proving the global existence of weak solutions of the regularized system \eqref{eq:poissonregeta}--\eqref{eq:gradientflowcontreg}
we shall need the following maximum principle for a semilinear PDE.


\begin{lemma}\label{lem:nonnegsubsol}
	Let $\Omega$ be an open, bounded subset of $\R^d$. For a fixed $T>0$ denote $\Omega_T := (0,T]\times \Omega$ and
	\begin{align*}
	C^2_1(\Omega_T) := \{ u\colon \Omega_T\to\R~ |~ u,\nabla u, \nabla^2 u, \partial_t u\in C(\Omega_T) \}.
	\end{align*}
	Let $\gamma>1$ and let $u\in C^2_1(\Omega_T)\cap C(\overline{\Omega}_T)$ be the classical solution of the initial/boundary-value problem
	\begin{align}\label{eq:initialboundaryvaluehelp}
	\begin{cases}
	\partial_t u=D^2\Delta u-\nu |u|^{\gamma-2}u &\text{in~}\Omega_T,\\
	\hphantom{\partial_t} u=0 &\text{on~}  [0,T]\times \partial\Omega,\\
	\hphantom{\partial_t}u=g &\text{on~} \{t=0\}\times \partial\Omega, 
	\end{cases}
	\end{align}
	with the nonnegative initial datum $g\colon \Omega \to \R$. Then, 
	\begin{align}\label{eq:maximumprinciple}
	\min_{\overline{\Omega}_T} u\geq 0.
	\end{align}
\end{lemma}

\begin{proof}
	Denote $U_T := \{(t,x)\in\Omega_T~|~u(t,x)<0\}$. Then $U_T$ is an open bounded subset of $\Omega_T$ and 
	\[
	   \partial_t u - D^2 \Delta u = -\nu |u|^{\gamma-2}u > 0\qquad \text{in~} U_T.
	\]
	Then using the classical weak maximum principle for the heat equation, see, e.g., \cite{Evans}, we have
	\[
	   \min_{\overline{U}_T} u = \min_{\partial{U}_T} u = 0.
	\]
	Consequently, $U_T = \emptyset$ and \eqref{eq:maximumprinciple} holds.
\end{proof}

\begin{lemma}
	Let $S\in L^2(\Omega)$ and $c^I\in H_0^1(\Omega)^{d\times d}\cap L^{\gamma}(\Omega)^{d\times d}$.
	Then for each $\eps>0$ the regularized system \eqref{eq:poissonregeta}--\eqref{eq:gradientflowcontreg}
	subject to the data \eqref{eq:contbc}--\eqref{eq:continitialdata} admits a global weak solution $(c,p)$
	satisfying \eqref{weaksol}. The regularized energy \eqref{eq:energycontregeta} satisfies
	\(  \label{regen}
	   \E^\eps[c(t)]+\sum_{k=1}^d \int_0^t \int_{\Omega} \bl \partial_t c^k(s,x)\br^2 \di x\di s = \E^\eps[c^I]\quad \text{for all~}t\geq 0.
	\)
\end{lemma}

\begin{proof}
We proceed along the lines of the proof of Theorem 2 of \cite{bionetworksanalysis}.
We employ the Leray-Schauder fixed point theorem in the space $L^2((0,T)\times\Omega)$.
For a given diagonal tensor $c\in L^2((0,T)\times\Omega)$ with nonnegative elements
we construct a solution $p_\eps\in H^1(\Omega)$ of the regularized Poisson equation \eqref{eq:poissonregeta} with no-flux boundary data
using the Lax-Milgram theorem; note that for $\eps>0$ the permeability tensor \eqref{Peps}
satisfies $\P^\eps \in L^\infty(\Omega)$, and uniform ellipticity follows from the assumption
$r\geq r_0>0$ in $\Omega$. Consequently, we have the uniform bound
\(  \label{peps_uniform}
   \Norm{\grad p^\eps}_{L^2(\Omega)} \leq C_\Omega \Norm{S}_{L^2(\Omega)}\qquad\mbox{for all }t\geq 0,\, \eps>0,
\)
where the constant $C_\Omega$ depends only on the domain $\Omega$; in particular, it is independent of $\eps>0$ and $c\in L^2((0,T)\times\Omega)$.

Existence of weak solutions $c_\eps$ of \eqref{eq:gradientflowcontreg} is obtained by a slight adaptation of Lemma 3 of \cite{bionetworksanalysis},
noting that for $\grad p \in L^2(\Omega)$ and $\eps>0$ the terms $\left( \partial_{x_k} p\right)^2\ast \eta_\eps$ are bounded in $L^\infty(\Omega)$.
The nonnegativity of the diagonal entries of $c_\eps$ follows from the fact that solutions of the semilinear PDE
\[
   \partial_t u=D^2\Delta u-\nu |u|^{\gamma-2}u
\]
are subsolutions to \eqref{eq:gradientflowcontreg}. Preservation of nonnegativity of $u$ for nonnegative initial and boundary data
has been established in Lemma \ref{lem:nonnegsubsol}.

The proof of continuity and compactness of the Schauder fixed point mapping $c \mapsto p_\eps \mapsto c_\eps$ in the space $L^2((0,T)\times\Omega)$
goes again along the lines of Theorem 2 of \cite{bionetworksanalysis}, using the so-called weak-strong lemma
for the Poisson equation (Lemma 7 of \cite{bionetworksanalysis}) and compact Sobolev embedding $H^1(\Omega) \subset L^2(\Omega)$.

The energy identity \eqref{regen} follows by multiplying the Poisson equation \eqref{eq:poissonregeta} by $p$ and integrating by parts,
\begin{align*}
	\sum_{k=1}^d \int_{\Omega} \bl r+c^k\ast\eta_\eps \br \bl \partial_{x_k} p\br^2 \di x=\int_{\Omega} S p \di x.
\end{align*}
Subtracting this from \eqref{eq:energycontregeta} we obtain
\begin{align*}
	\E^\eps[c]=\sum_{k=1}^d \bl \int_{\Omega} \frac{D^2}{2}  \left|\grad c^k\right|^2 - \bl r+c^k\ast\eta_\eps \br  \bl\partial_{x_k} p\br^2 
	  + \frac{\nu}{\gamma} \bigl| c^k \bigr|^\gamma  \d x\br+2\int_{\Omega} S p\di x.
\end{align*}
Integration by parts in suitable terms and using \eqref{eq:poissonregeta} yields then
\begin{align*}
  \frac{\di}{\di t} \E^\eps[c]&=\sum_{k=1}^d \bl \int_{\Omega} -D^2 \Delta c^k \partial_t c^k+2\partial_{x_k} \bl \bl r+c^k\ast\eta_\eps \br \partial_{x_k} p\br \partial_t p
     -\partial_t c^k\bl\partial_{x_k} p\br^2 +\nu  \bigl| c^k \bigr|^{\gamma-1}\partial_t c^k\di x\br \\
     & \qquad+2\int_{\Omega} S \partial_t p\di x\\
     & = -\sum_{k=1}^d \int_{\Omega} \bl D^2 \Delta c^k +\bl\partial_{x_k} p\br^2\ast\eta_\eps  -\nu \bigl| c^k \bigr|^{\gamma-1}\br \partial_t c^k\di x
       +2\int_{\Omega}\bl \nabla\cdot \bl \mathbb{P}^\eps[c] \nabla p \br +S \br\partial_t p\di x\\&= -\sum_{k=1}^d \int_{\Omega} \bl \partial_t c^k\br^2 \di x.
\end{align*}
and an integration in time gives \eqref{regen}.
\end{proof}

\noindent
The passage to the limit $\epsilon\to 0$ in \eqref{eq:poissonregeta}--\eqref{eq:gradientflowcontreg}
is based on the uniform apriori estimates
\[
      c\in L^{\infty}(0,\infty;H_0^1(\Omega))\cap L^{\infty}(0,\infty;L^{\gamma}(\Omega)),\quad \partial_t c\in L^2((0,\infty)\times\Omega),\\
      \grad p\in L^{\infty}(0,\infty; L^2(\Omega)),\quad \sqrt{c^k\ast\eta_\eps}\partial_{x_k} p\in L^{\infty}(0,\infty;L^2(\Omega)),\; k=1,\dots,d,
\]
which follow from the energy identity \eqref{regen} and from \eqref{peps_uniform}.
Then, since a subsequence of $c^\eps\ast\eta_\eps$ converges strongly to $c$ in the norm topology of $L^2((0,T)\times \Omega)$,
a slight modification of Lemma 7 in \cite{bionetworksanalysis} gives the strong convergence of $p^{\eps}$ to $p$ in $L^2(0,T;H^1(\Omega))$ with no-flux boundary data
where $p$ is the unique solution of the Poisson equation \eqref{PcontR2} with given $c$.
Thus, $\bl \partial_{x_k} p^\eps\br^2$ converges strongly  to $\bl \partial_{x_k} p\br^2$ in $L^1\bl (0,T)\times \Omega\br$
and $\bl\partial_{x_k} p^\eps\br^2 \ast\eta_\eps$ also converges strongly to $\bl \partial_{x_k} p\br^2$ in $L^1\bl (0,T)\times \Omega\br$.
The limit passage in the metabolic term $\bigl| c^k\bigr|^{\gamma-2} c^k$ can be shown as in Lemma 4 in \cite{bionetworksanalysis}
due to the uniform boundedness of $c^{\eps}$ in $L^{\gamma}((0,T)\times \Omega)$.
The energy dissipation inequality \eqref{eq:energydissinequ} follows by passing to the limit $\eps\to 0$
in \eqref{regen} using the weak lower semicontinuity of the $L^2$-norm.
This concludes the proof of Theorem \ref{th:existencecont}.

\section{Numerical simulations}\label{sec:numerics}
In this section we provide results of numerical simulations for the discrete model introduced in Section \ref{sec:micro}.
We implement a minimization scheme for the discrete energy \eqref{eq:energydisc} constrained by the Kirchhoff law \eqref{eq:kirchhoff},
based on the numerical methods proposed in \cite{bionetworksbookchapter}.

\begin{figure}[htbp]
	\centering
	\subfloat[Initial data] {\includegraphics[width=0.24\textwidth]{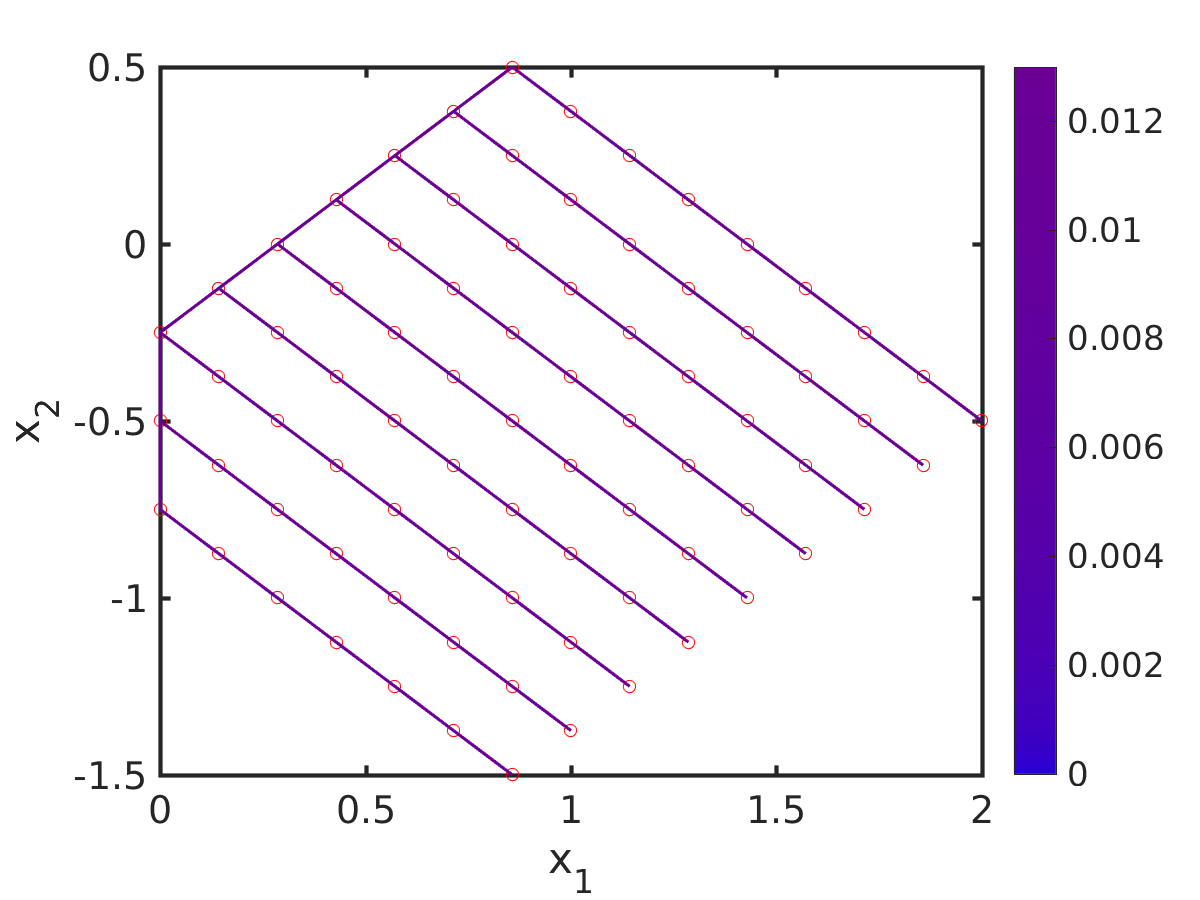}\label{fig:initialcond}}
	\subfloat[$\epsilon=10^{-4}$] {\includegraphics[width=0.24\textwidth]{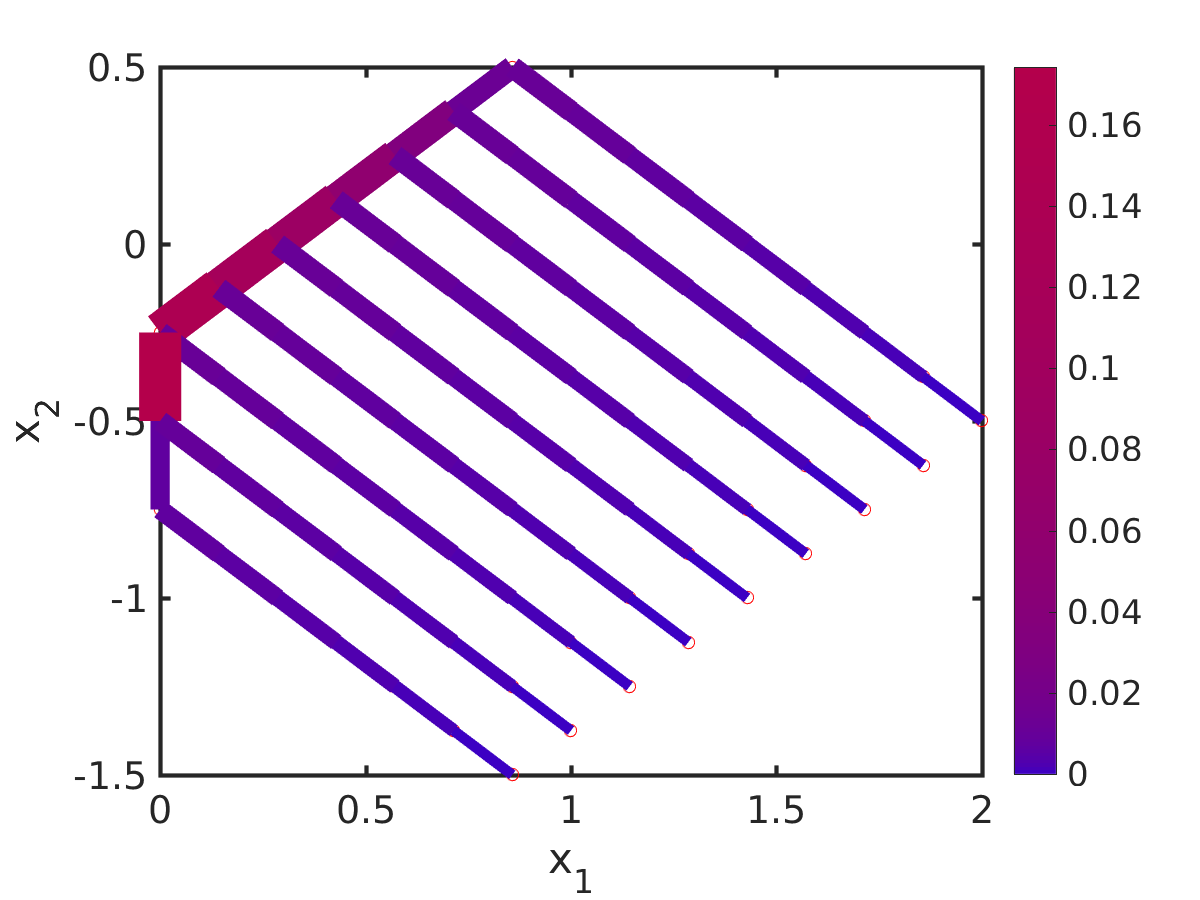}\label{fig:perturb-4}}
	\subfloat[$\epsilon=10^{-3}$] {\includegraphics[width=0.24\textwidth]{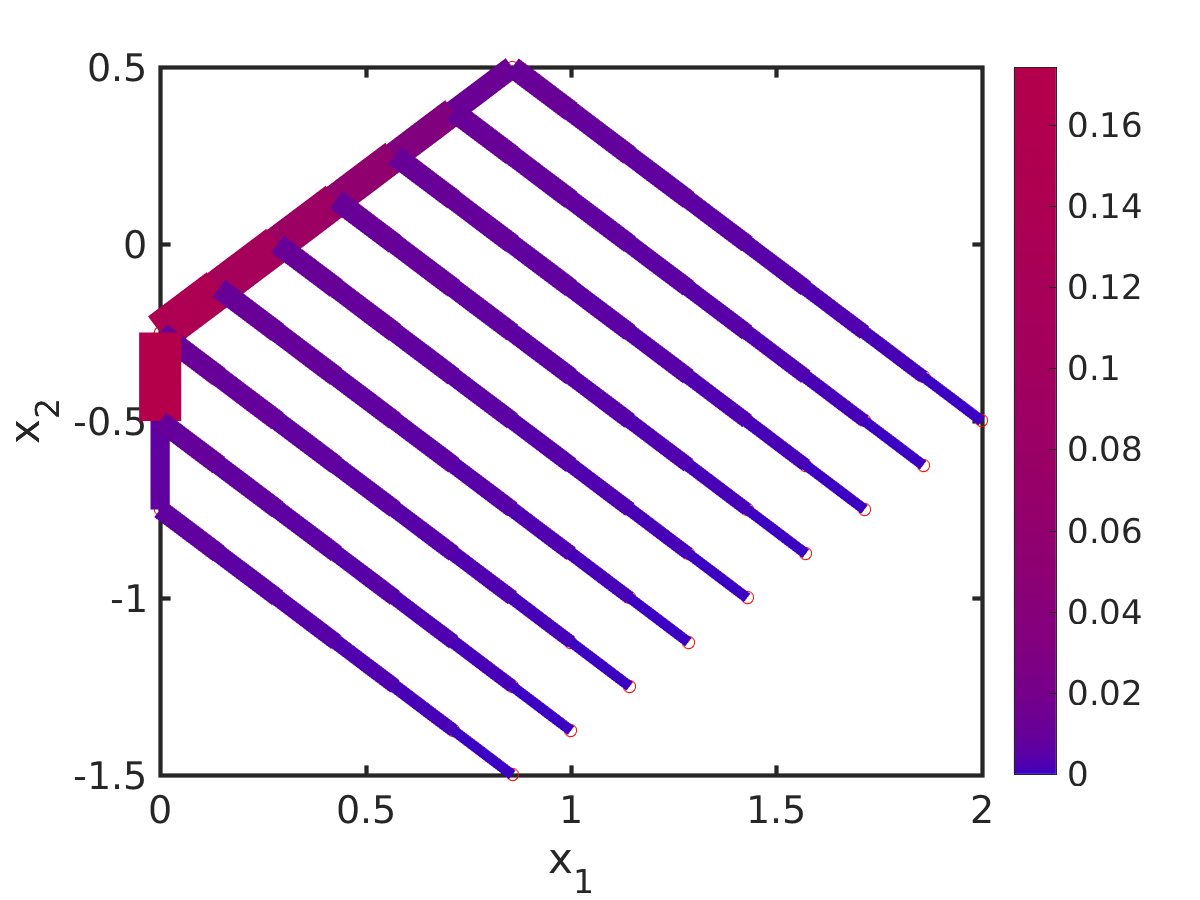}}
	\subfloat[$\epsilon=10^{-2}$] {\includegraphics[width=0.24\textwidth]{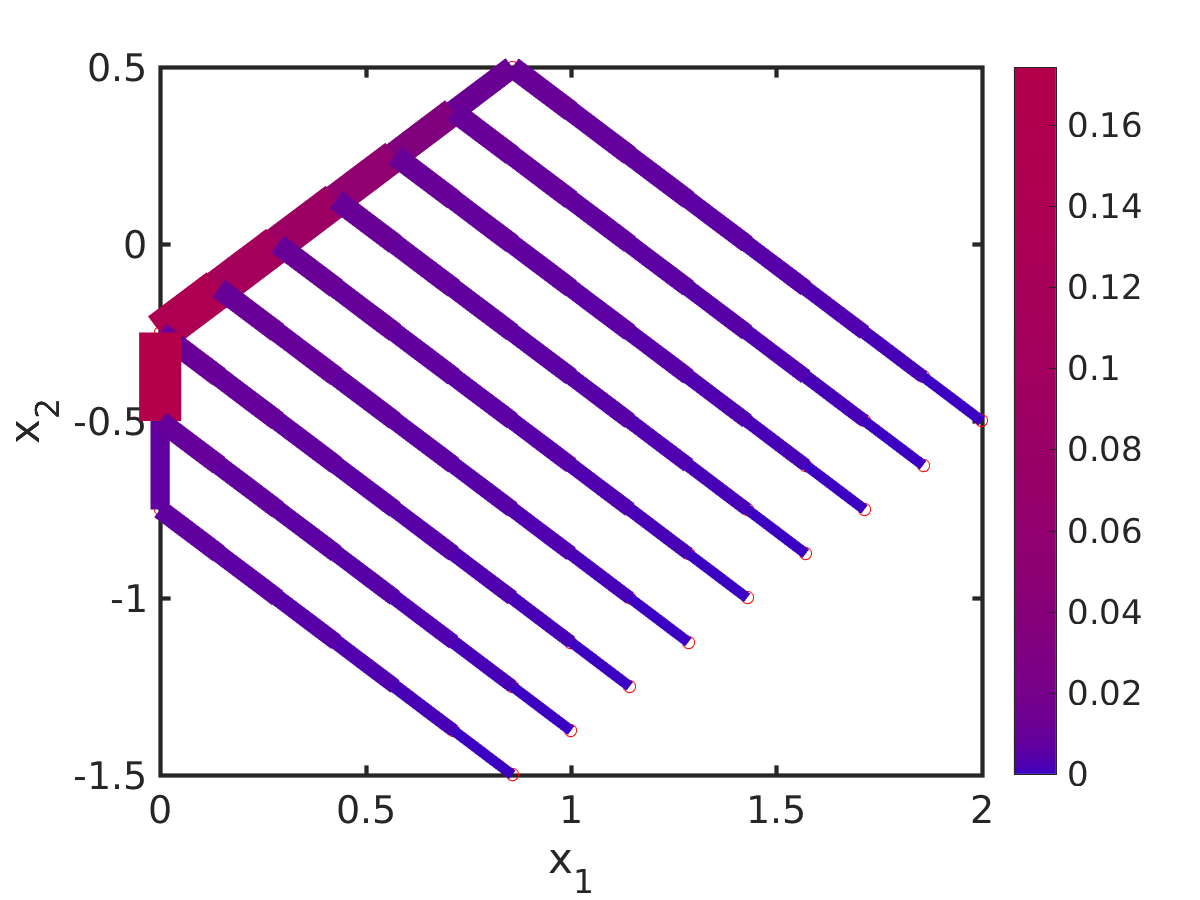}}\\
	\subfloat[$\epsilon=0.1$] {\includegraphics[width=0.24\textwidth]{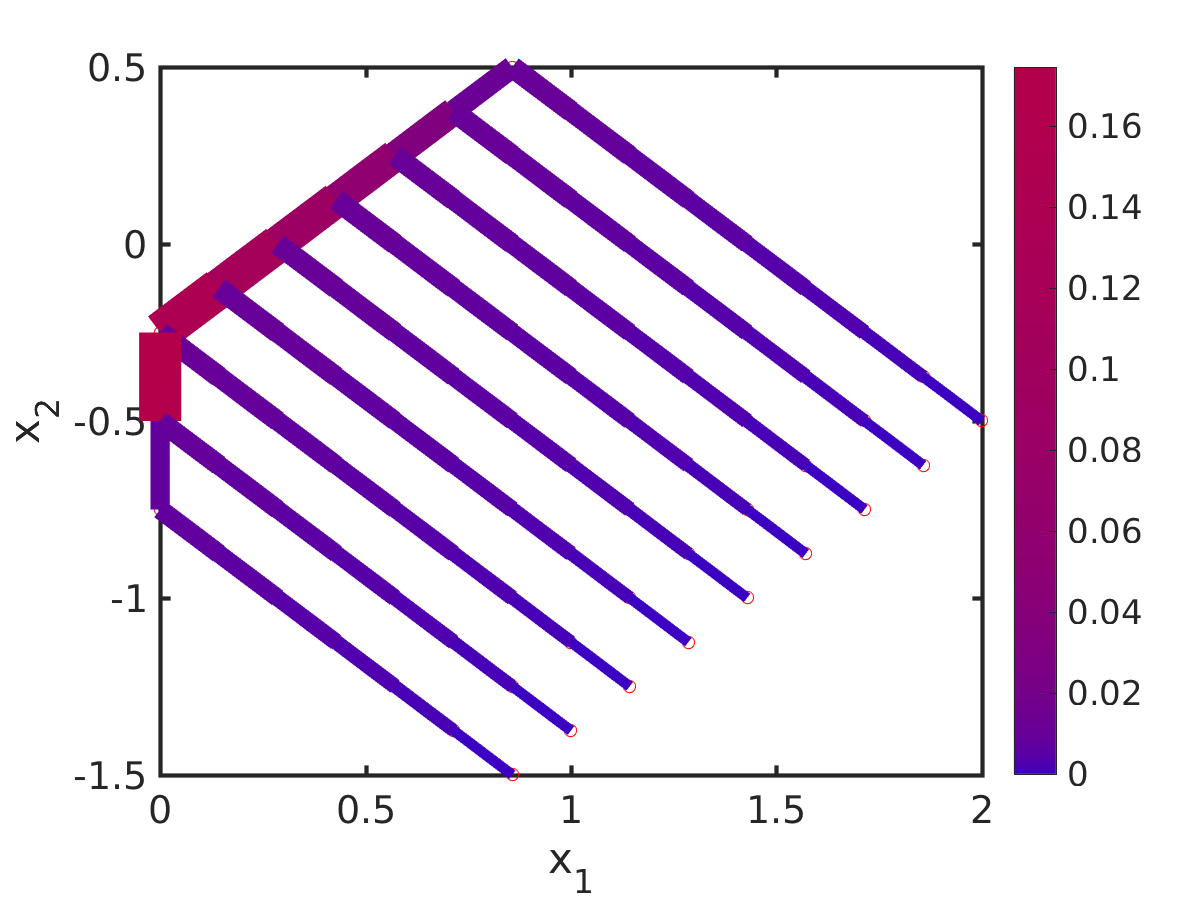}}
	\subfloat[$\epsilon=0.5$] {\includegraphics[width=0.24\textwidth]{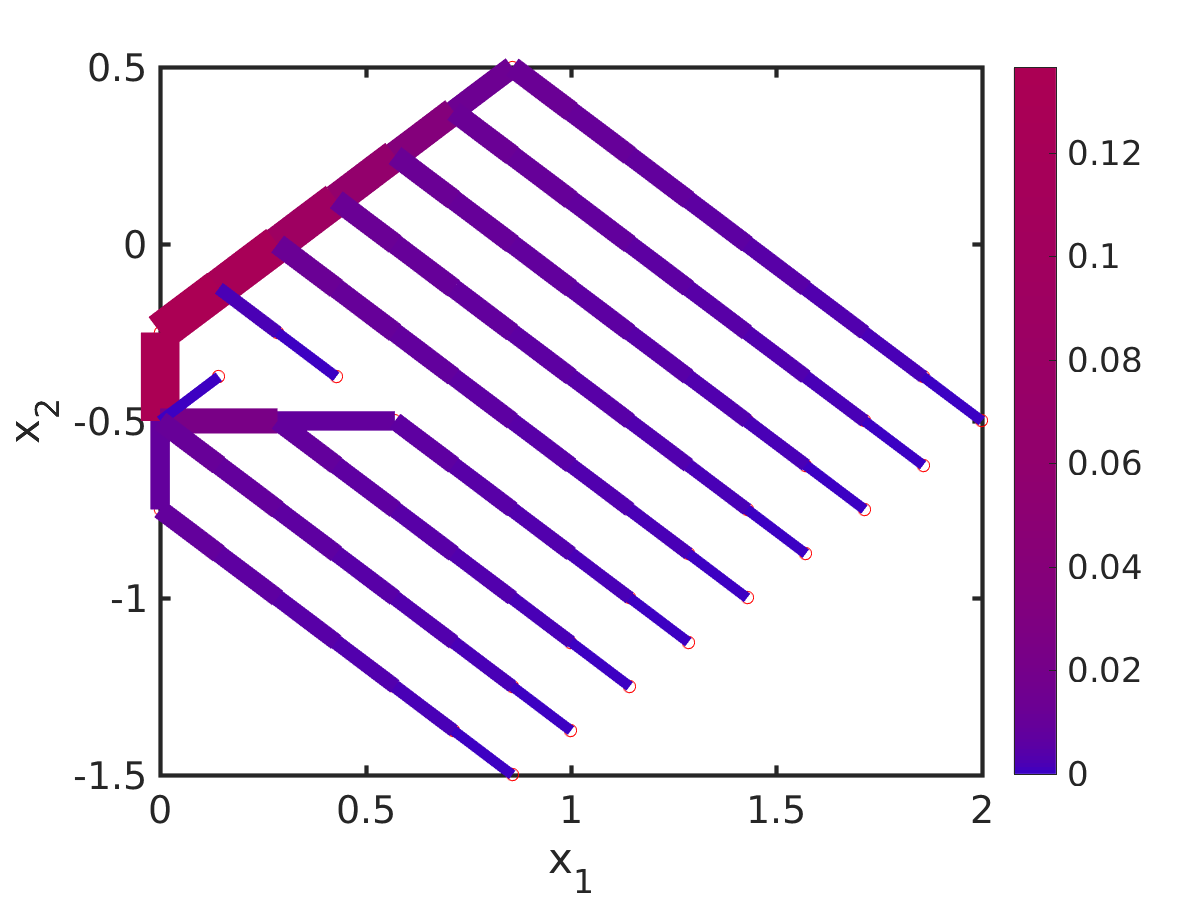}}
	\subfloat[$\epsilon=1$] {\includegraphics[width=0.24\textwidth]{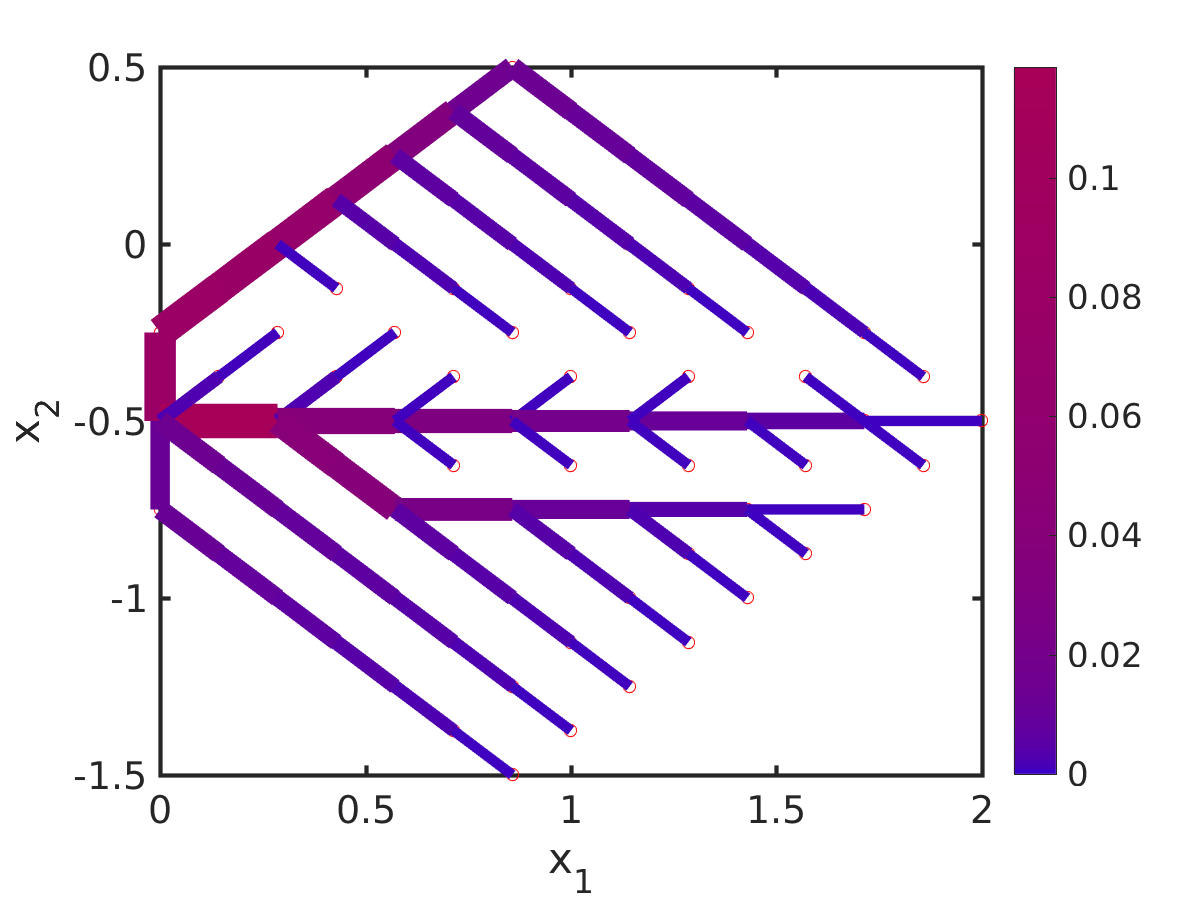}}
	\subfloat[$\epsilon=2$] {\includegraphics[width=0.24\textwidth]{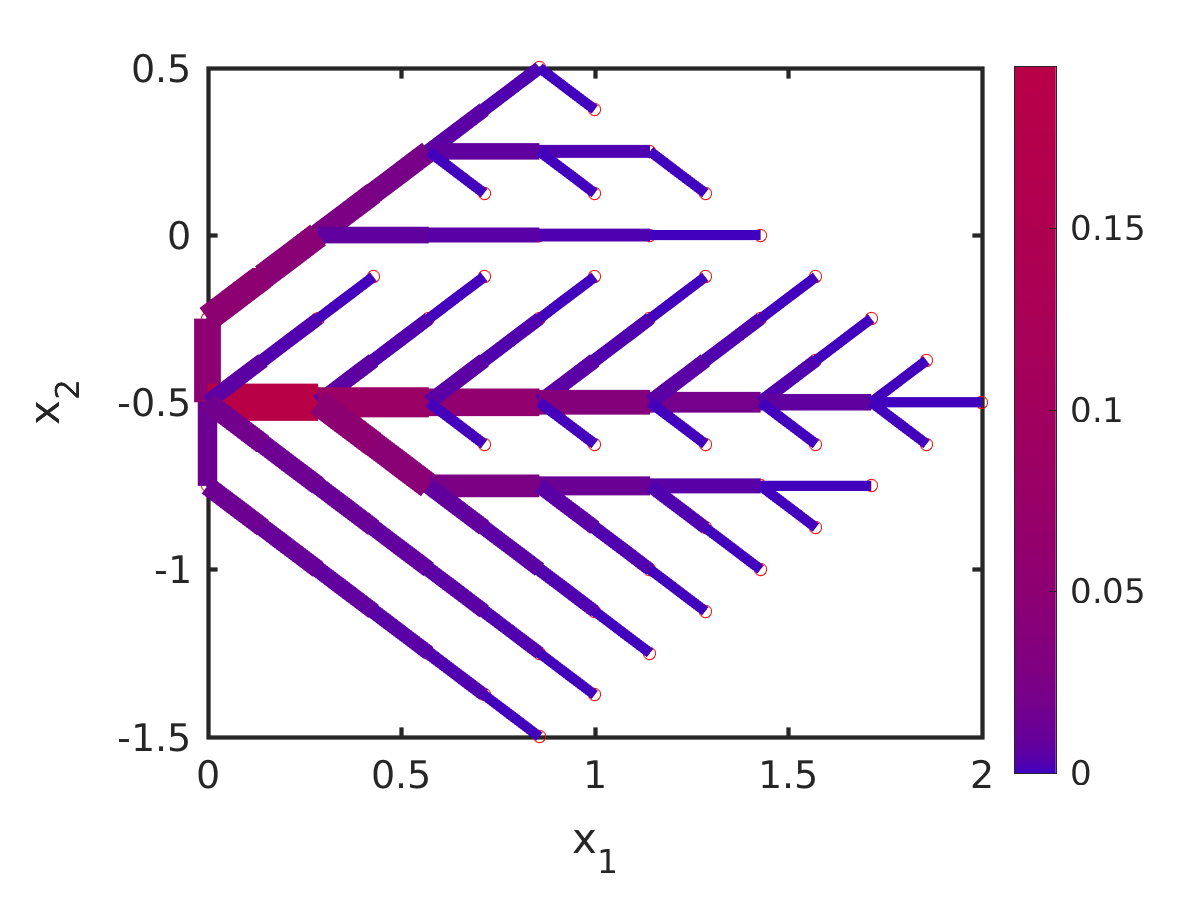}}
	\caption{Stability of steady states under perturbations $\epsilon$ of initial data for the discrete model.}\label{fig:perturbationinitialcond}
\end{figure}

For the numerical simulations we consider a planar graph $G=(\Vset,\Eset)$ whose vertices and edges define
a diamond shaped geometry embedded in the two-dimensional domain $\Omega=(0,2)\times(-1.5,0.5)$.
We consider $|\Vset|=78$ vertices and $|\Eset|=201$ edges.
For vertex $i\in\Vset$ let $(x^i,y^i)$ denote its position.
The source $S$ is assumed to be positive on the subset of vertices
$$\Vset^+:=\{i\in\Vset;\; x^i\leq 0.1\}$$
and constant and negative on its complement $\Vset\backslash\Vset^+$. For $i\in\Vset$ we set
\begin{align*}
   S_i:=\begin{cases}
       \sigma_i^+ & i\in \Vset^+\\
       \sigma_i^- & i\in \Vset\backslash\Vset^+
   \end{cases}
\end{align*}
where
\begin{align*}
\sigma_i^+ := 10^4 \exp\bl -10\bl 50x_i^2+10\bl y_i+0.5\br^4\br\br,\quad \sigma_i^- := -\frac{1}{|\Vset\backslash \Vset^+|}\sum_{j\in\Vset^+}\sigma_j^+.
\end{align*}
In the sequel we prescribe the initial condition $\overline{C}=\bl \overline{C}_{ij}\br_{(i,j)\in\Eset}$, unless stated otherwise. We assume
$\overline{C}_{ij}:=5$ for every $(i,j)\in\Eset$ on a tree, see Figure \subref*{fig:initialcond},
and $\overline{C}_{ij}:=10^{-10}$ otherwise.

For solving the constrained energy minimization problem we consider the following iterative procedure:
\begin{itemize}
\item Initialization: For each edge $(i,j)\in \Eset$ compute its length $L_{ij}$ 
and define the parameters $\nu:=1$, $\tau:=0.025$ and $tol:=10^{-6}$.

\item Step 1 (Pressure): For $\overline{C}$ given, compute the coefficient matrix $B=(b_{ij})\in\R^{n-1,n-1}$ with entries
\begin{align}\label{eq:coefficientmatrixkirchhoff}
b_{ij}&=\begin{cases}
-\frac{C_{ij}}{L_{ij}^2} & \quad(i,j)\in\Eset\\
0 & \quad(i,j)\notin\Eset
\end{cases},\quad i,j=1,\ldots,n-1,\quad i\neq j,\\
b_{ii}&=\sum_{j\in N(i)}\frac{C_{ij}}{L_{ij}^2},\quad
i=1,\ldots,n-1.
\end{align}
and solve via least square minimization:
\begin{align*}
\min_P\|BP-S\|_2 
\end{align*}

\item Step 2 (Conductivity): For given pressure $P$ and conductivities $\overline{C}$ find a minimizer $C$ of the regularization
\begin{align}\label{eq:energydiscnumericsreg}
   \E^{\tau}[C] := \frac{\|C-\overline{C}\|_2^2}{2\tau} + \sum_{(i,j)\in\Eset}\bl \frac{Q_{ij}(C)^2}{C_{ij}}+\nu C_{ij}^{\gamma}\br L_{ij}
\end{align}
 of the discrete energy functional \eqref{eq:energydisc} via interior point method for a regularisation parameter  $\tau>0$.
  
\item Step 3 (Energy decrease): If $\left| \E^{\tau}[C]-\E^{\tau}[\overline{C}]\right|>tol$, set $\overline{C}:=C$ and go back to step 1.
\end{itemize}

\begin{figure}[htbp]
	\centering
	\subfloat[$\gamma=0.5$] {\includegraphics[width=0.24\textwidth]{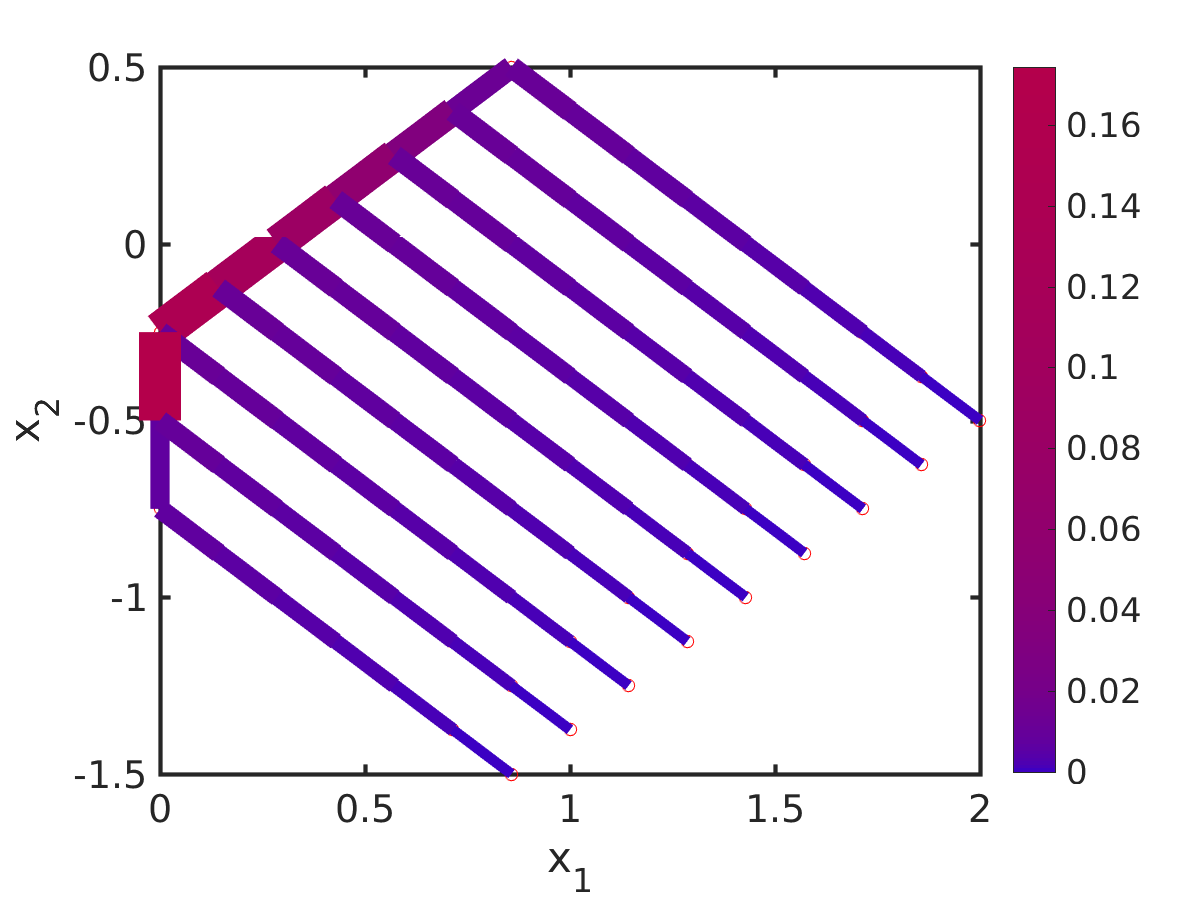}
\includegraphics[width=0.24\textwidth]{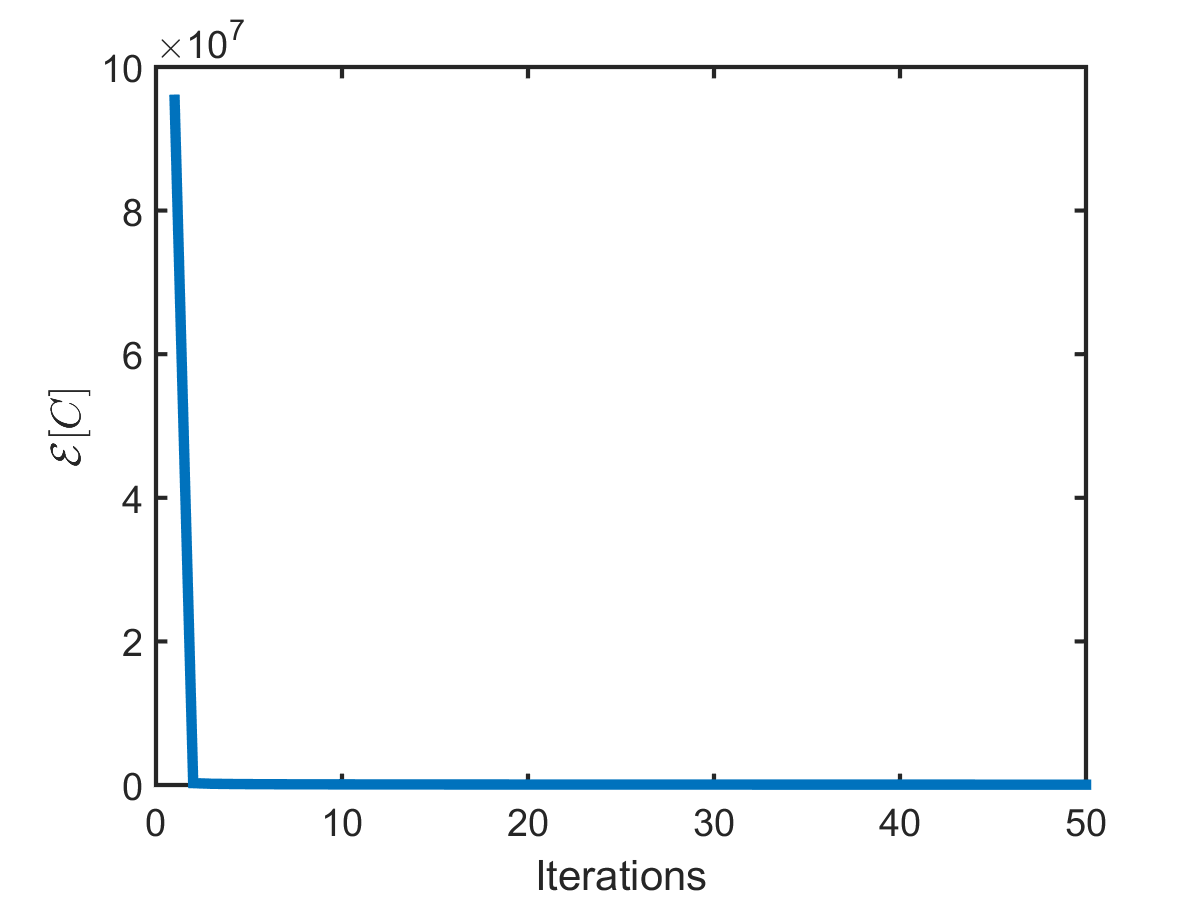}}
	\subfloat[$\gamma=1.5$] {\includegraphics[width=0.24\textwidth]{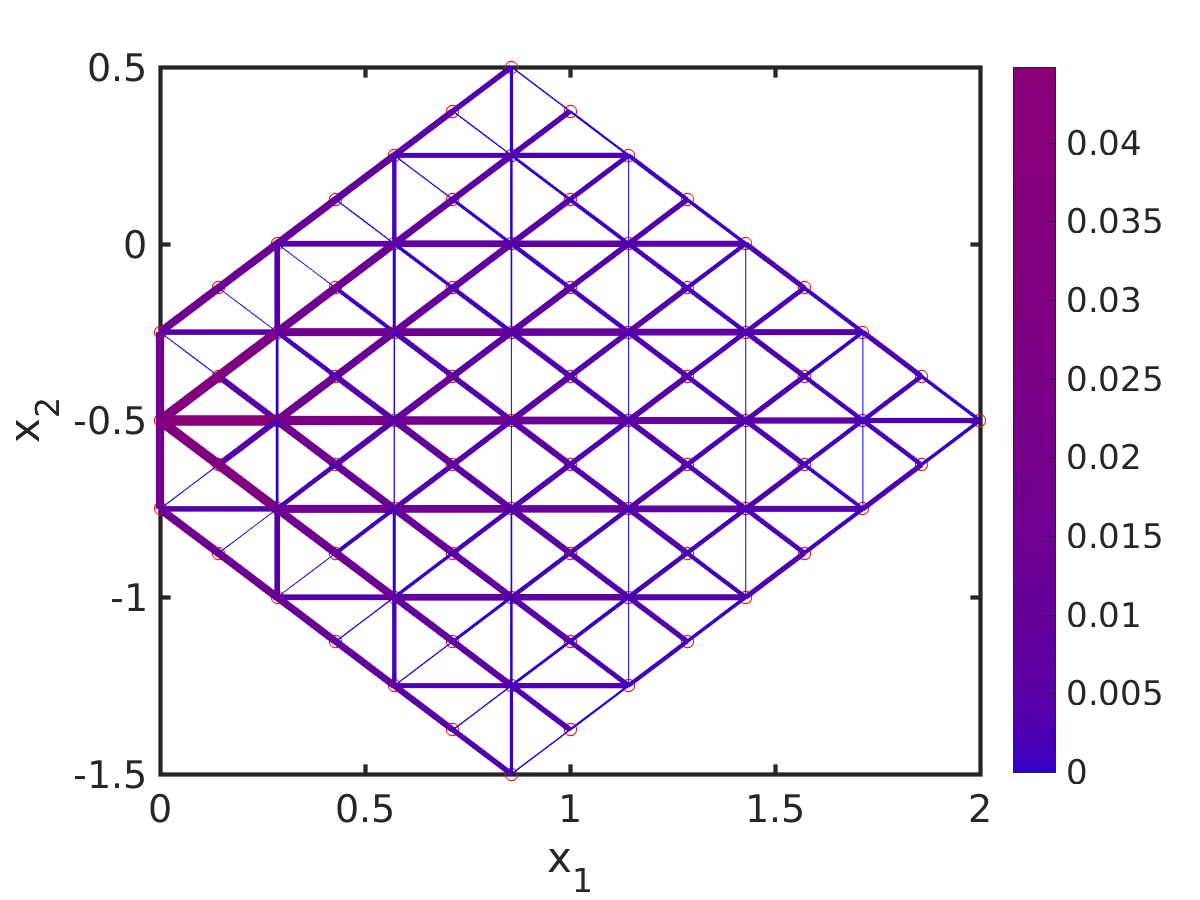}\includegraphics[width=0.24\textwidth]{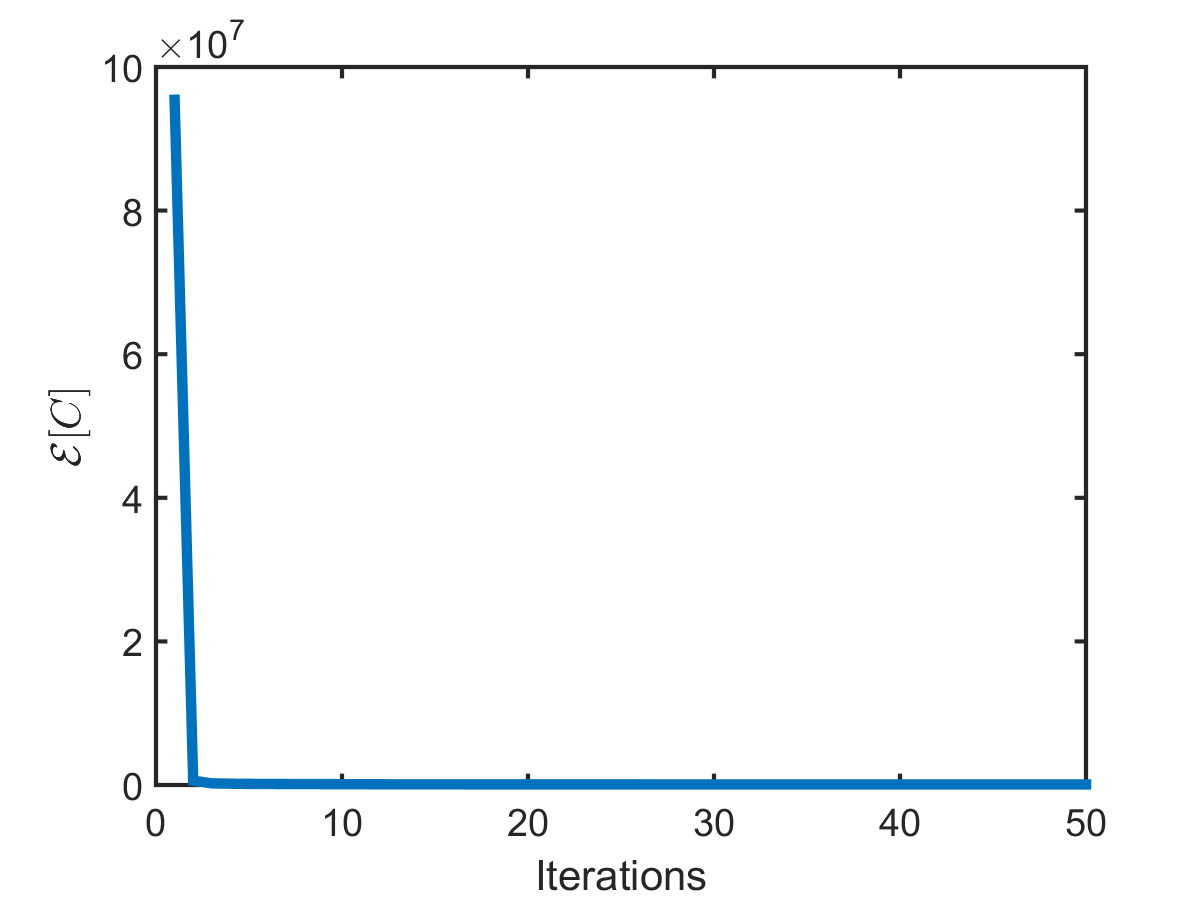}}	
	\caption{Stationary solution to the discrete model  and its  decrease of energy for different values of $\gamma>0$.}\label{fig:variationgamma}
\end{figure}

Note that for $\tau>0$ solving \eqref{eq:energydiscnumericsreg} is equivalent to an implicit Euler step for \eqref{eq:GF-alpha}.
The choice of the time step $\tau>0$ is crucial. On the one hand, the time step should not be chosen too large so that an accurate solution can be obtained.
On the other hand, choosing $\tau$ too small may result in very long simulation times, especially because the convergence seems to be very slow close to the minimizer,
compare Figure \ref{fig:variationgamma} where the slow decay of the energy functional is shown. Armijo's condition \cite{NoceWrig06} 
suggests a good choice of the parameter $\tau$ so that sufficient decrease of the energy functional is achieved in every time step.

In the sequel we present the energy minima (stationary solutions) obtained by the above algorithm for different values of $\gamma$.
For every edge $(i,j)\in\Eset$ we plot the value of the conductivity $C_{ij}$ in terms of the width of the associated edge.
In Figure \ref{fig:perturbationinitialcond} we show the steady states 
under an $\eps$-perturbation of the initial condition $\overline{C}$ for $\gamma=0.5$, i.e., we consider $\overline{C}_{ij}+\epsilon$
instead of $\overline{C}_{ij}$ for all edges $(i,j)\in\Eset$. As shown in Figure \ref{fig:perturbationinitialcond} the steady states are
the same trees for small perturbations, e.g., $\epsilon\leq 0.1$, as the tree given by the initial condition in Figure \subref*{fig:initialcond}.
In particular, the steady states are stable under small perturbations of the initial condition.
For larger perturbations, e.g., $\epsilon\in\{0.5,1,2\}$, we obtain  steady states different from the initial condition.
This illustrates that the energy functional \eqref{eq:energydisc} has multiple local minima
and, consequently, the system \eqref{eq:GF-alpha}--\eqref{eq:kirchhoff} has non-unique steady states. In particular, the steady states strongly depend on the choice of the initial data.

In Figure \ref{fig:variationgamma} the stationary solution of \eqref{eq:GF-alpha}--\eqref{eq:kirchhoff}
and the decay of the energy functional are shown for different values of $\gamma>0$.
Note that the stationary solution is a tree for $\gamma=0.5$ and a full network for $\gamma=1.5$. 
This is in agreement with the observations of \cite{hu} where a phase transition at $\gamma=1$ was suggested with steady states in the form of a tree for $\gamma<1$ and full networks as steady states for $\gamma>1$.

\begin{figure}[htbp]
	\centering
	\subfloat[Initial data] {\includegraphics[width=0.24\textwidth]{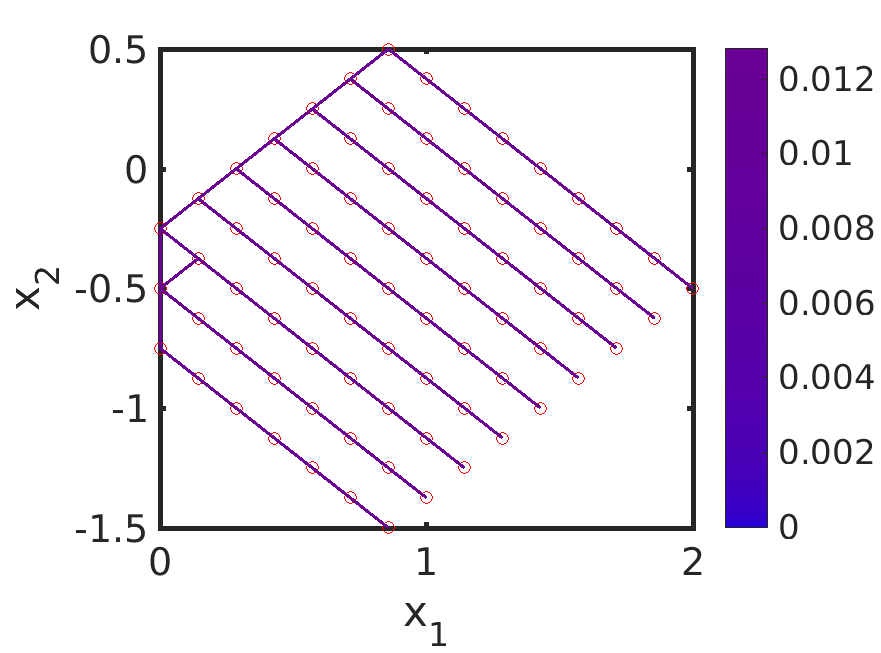}\includegraphics[width=0.24\textwidth]{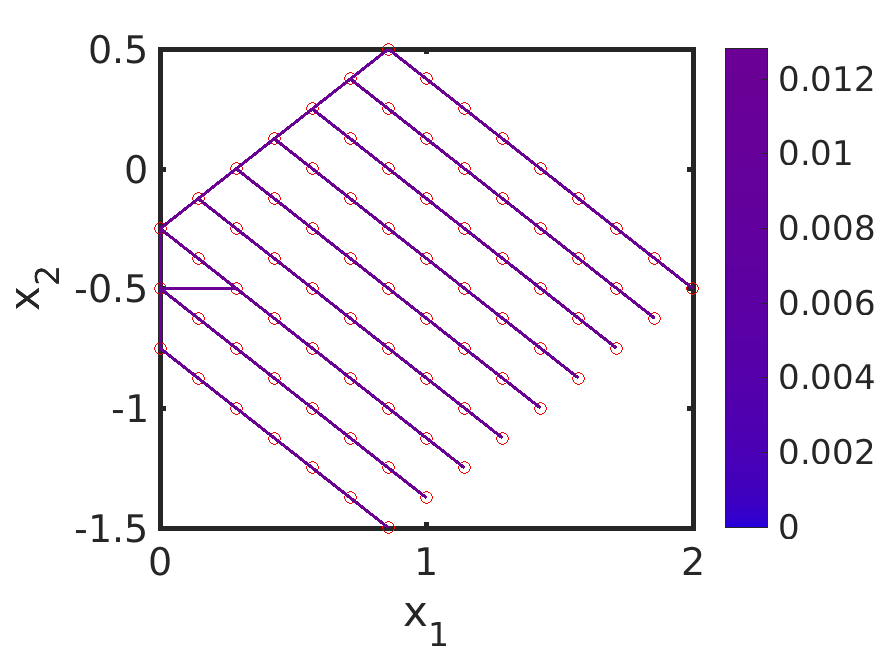}\includegraphics[width=0.24\textwidth]{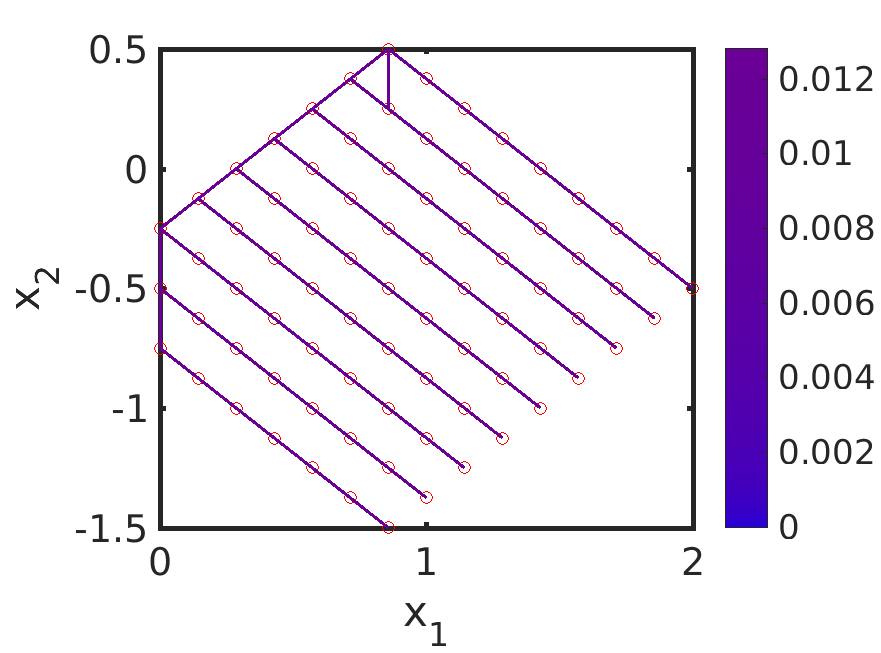}\includegraphics[width=0.24\textwidth]{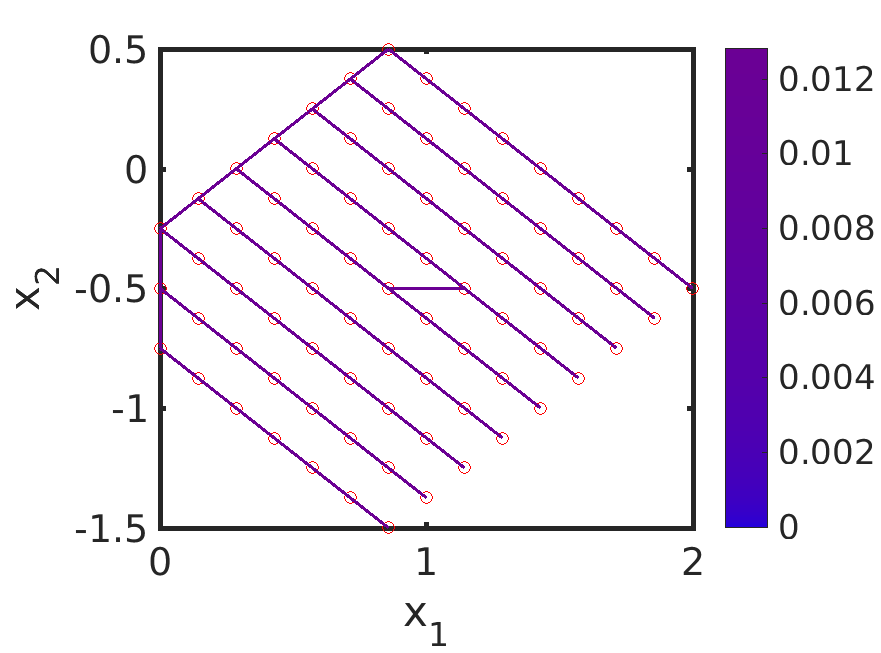}\label{fig:closeoneloopinitialcondinit}}
	
	\subfloat[Associated steady states] {\includegraphics[width=0.24\textwidth]{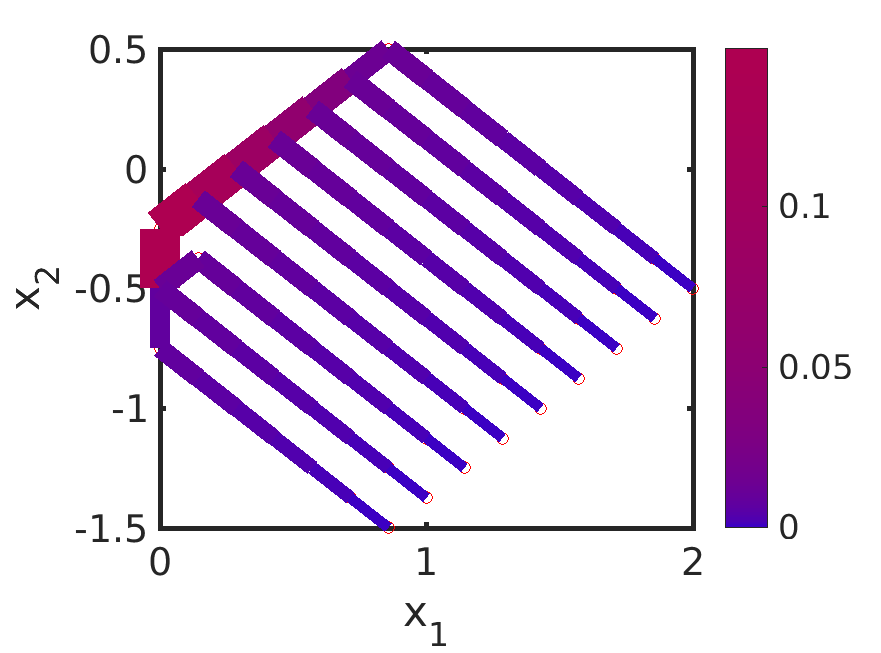}\includegraphics[width=0.24\textwidth]{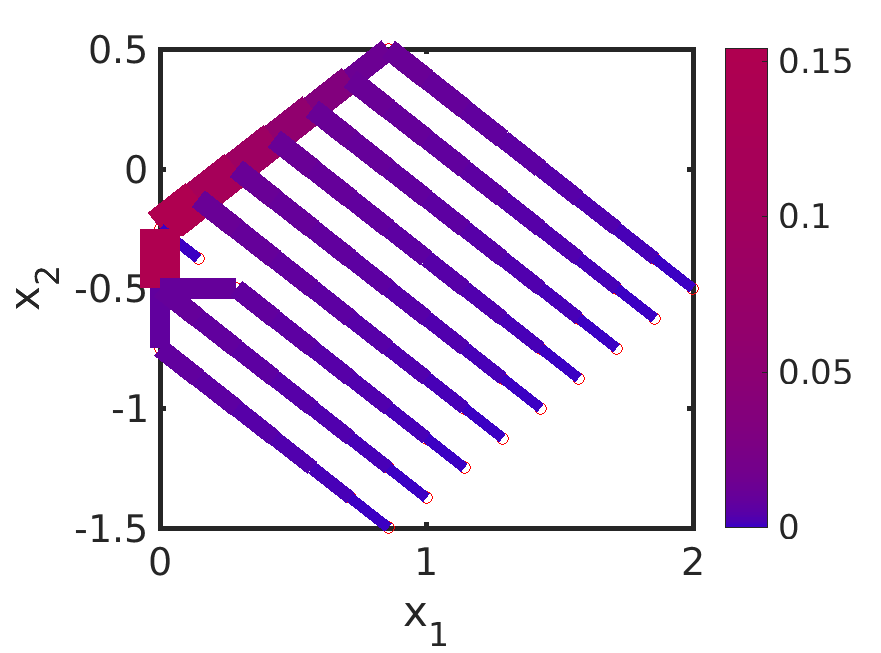}\includegraphics[width=0.24\textwidth]{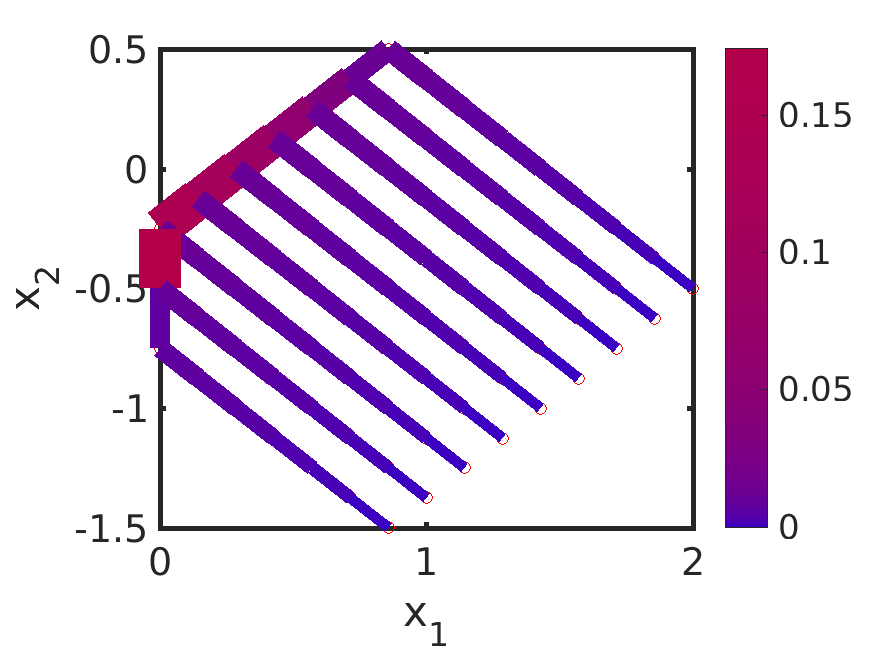}\includegraphics[width=0.24\textwidth]{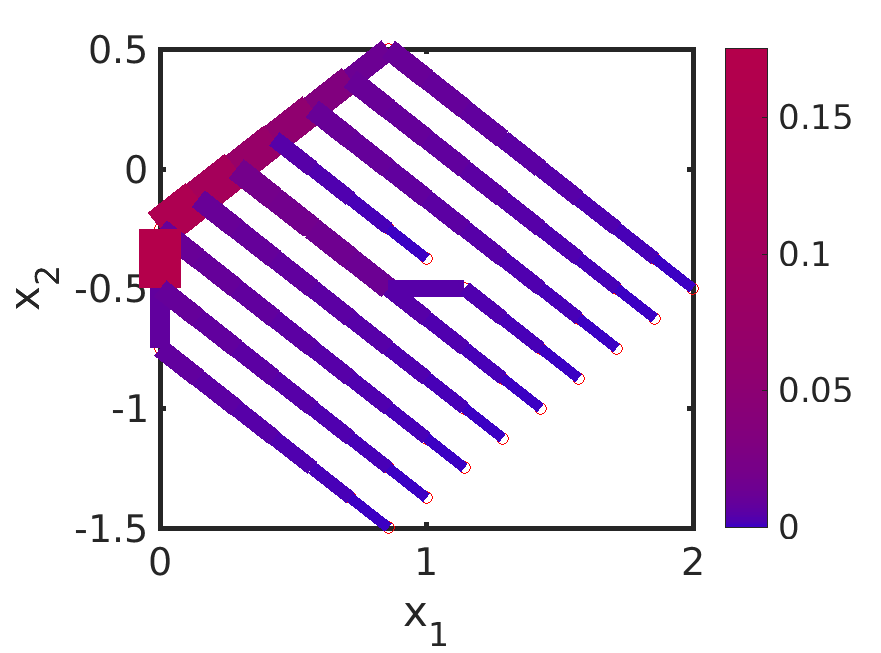}\label{fig:closeoneloopinitialcondsteadystate}}
	\caption{Stability of steady states when one loop in tree-structured initial data is closed.}\label{fig:closeoneloopinitialcond}
\end{figure}

In Figure \ref{fig:closeoneloopinitialcond} we consider initial data in form of a tree, Figure \subref*{fig:initialcond}, and close one of its loops, as shown in Figure \subref*{fig:closeoneloopinitialcondinit}. These initial conditions lead to the steady states in Figure \subref*{fig:closeoneloopinitialcondsteadystate}.
Note that closing one loop in the initial data 
leads to steady states which  only differ locally (i.e., in a neighborhood of the loop) from the original tree in Figure \subref*{fig:initialcond}. Closing one loop in areas of smaller conductivities in the associated steady state leads to the same tree structure as in the original initial data in Figure \subref*{fig:initialcond} as shown for the third choice of initial data in Figure \subref*{fig:closeoneloopinitialcondinit}. In particular, closing loops at different locations leads to different steady states in general, unless the resulting steady state is the tree in the initial condition in Figure \subref*{fig:initialcond}. This shows again that we obtain trees as steady states for $\gamma=0.5$, the steady states are non-unique and the form of the steady states strongly depends on the given initial data. In particular, loops in the initial data are opened over time for $\gamma=0.5$.

\begin{figure}[htbp]
	\centering
	\subfloat[Initial data] {\includegraphics[width=0.24\textwidth]{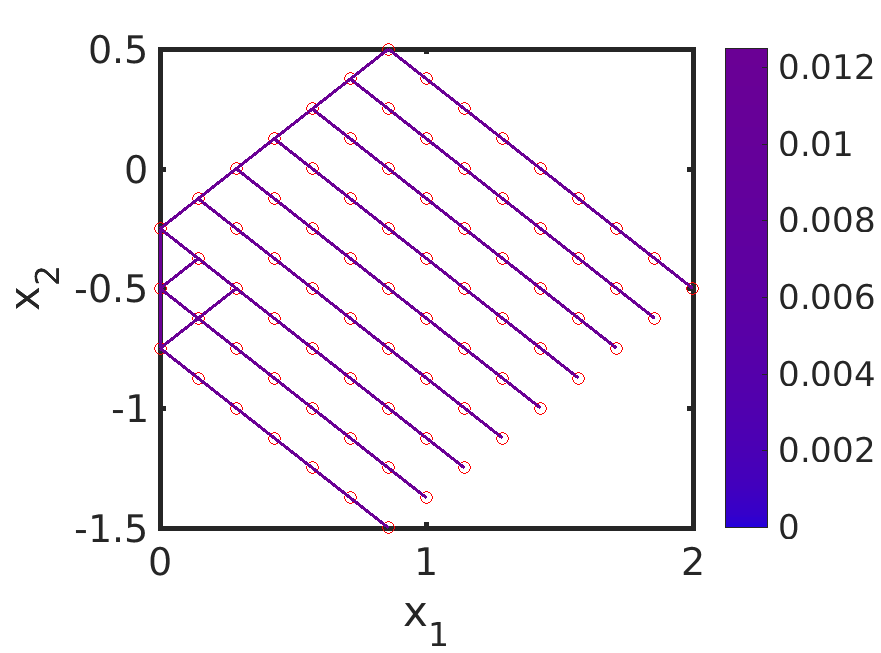}\includegraphics[width=0.24\textwidth]{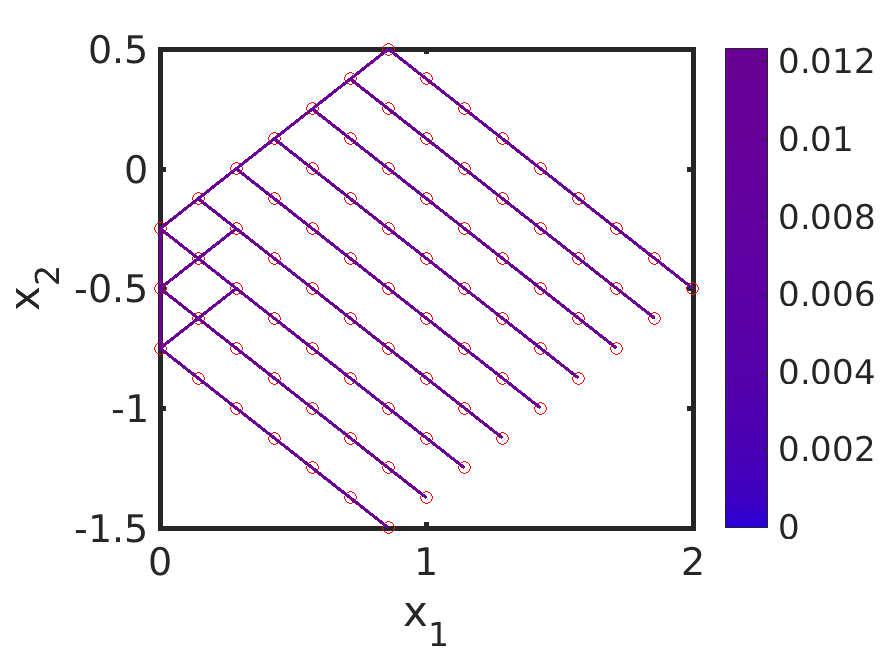}\includegraphics[width=0.24\textwidth]{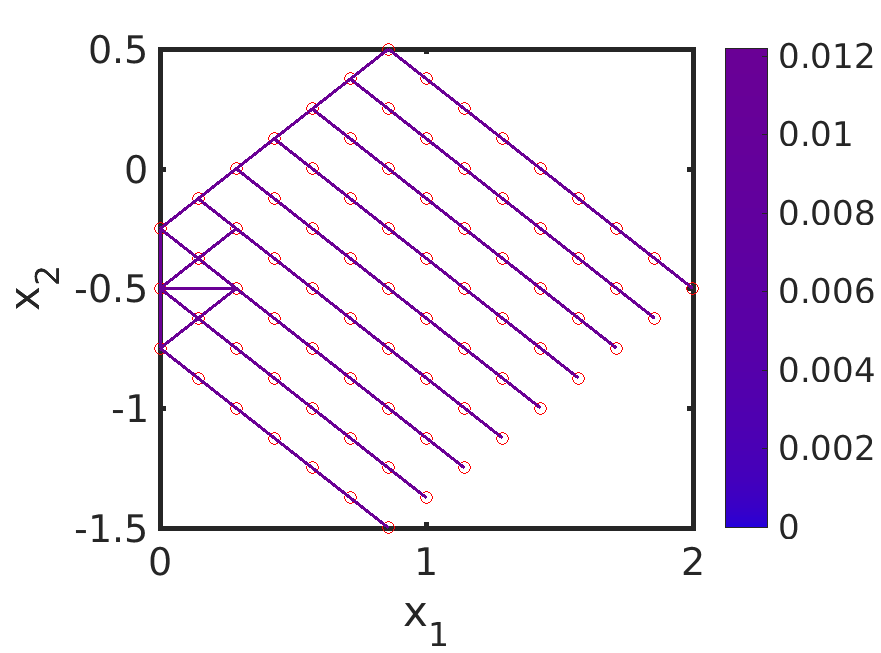}\includegraphics[width=0.24\textwidth]{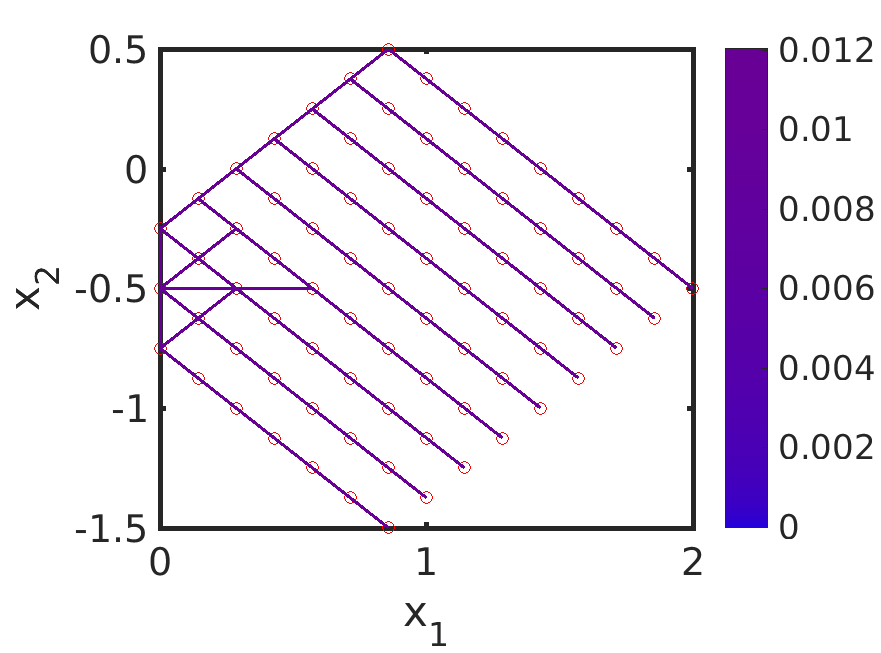}\label{fig:closeloopsleftinitialcondinit}}
	
	\subfloat[Associated steady states] {\includegraphics[width=0.24\textwidth]{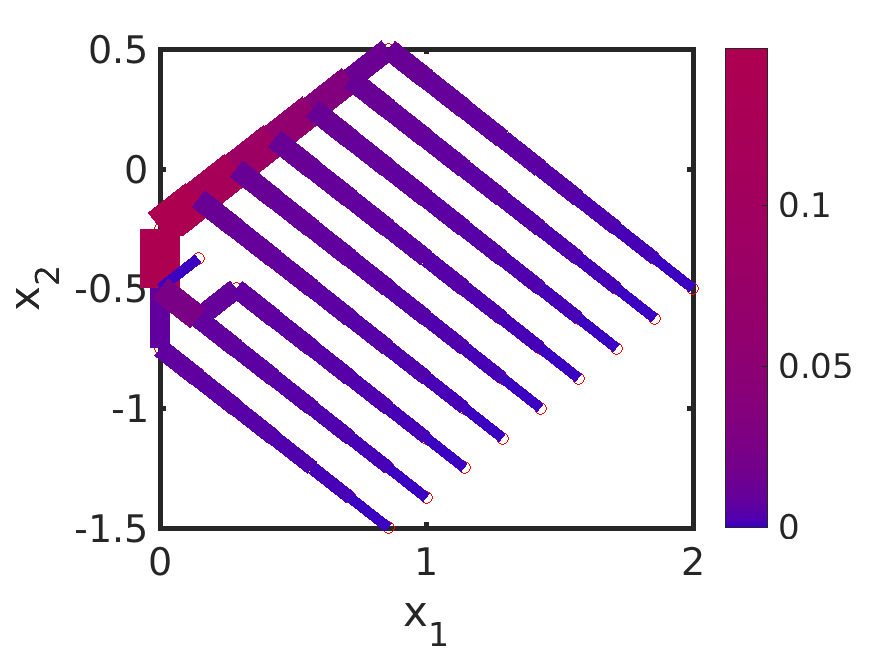}\includegraphics[width=0.24\textwidth]{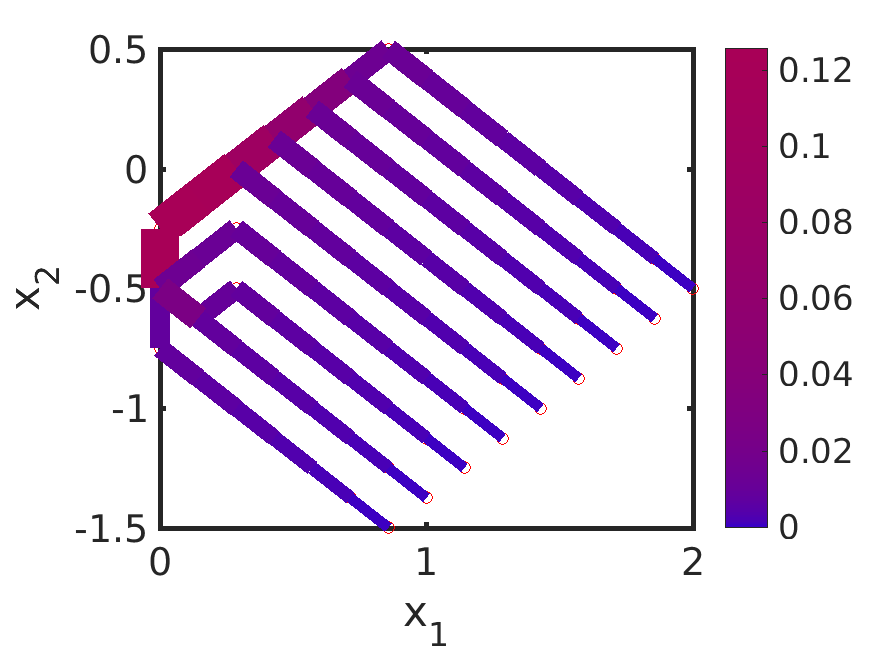}\includegraphics[width=0.24\textwidth]{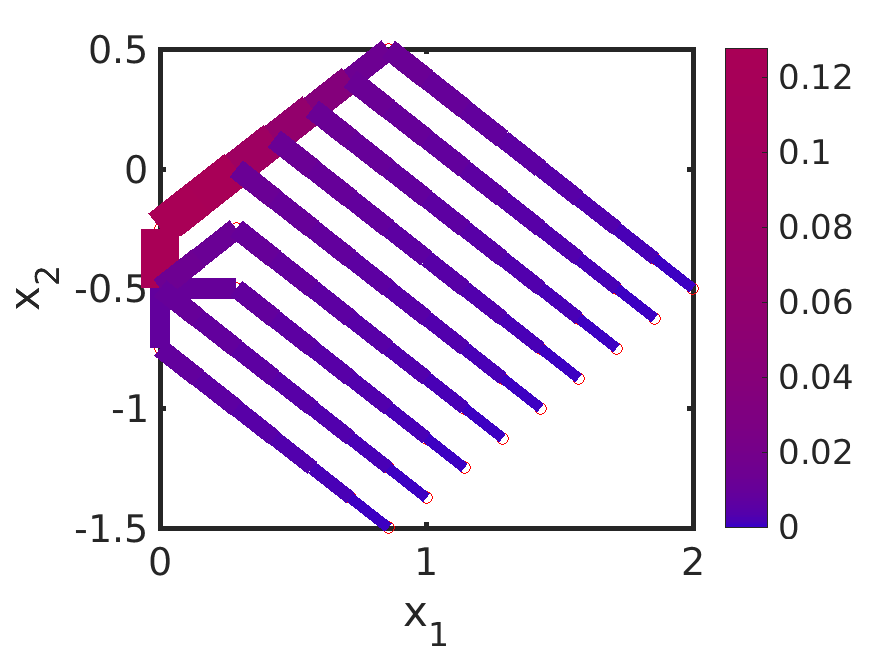}\includegraphics[width=0.24\textwidth]{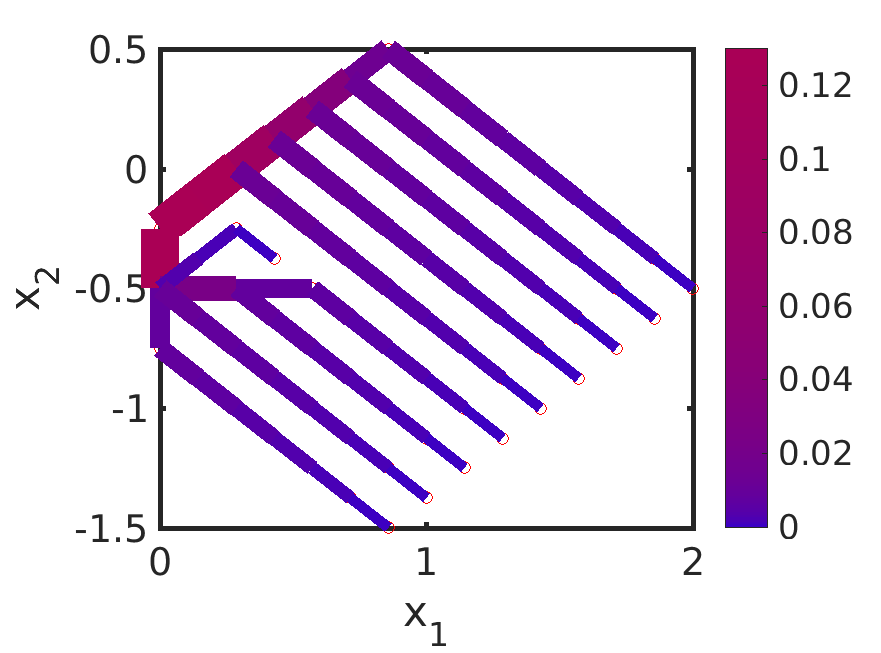}\label{fig:closeloopsleftinitialcondsteadystate}}
	\caption{Stability of steady states when several loops in tree-structured initial data are closed in the discrete model.}\label{fig:closeloopsleftinitialcond}
\end{figure}

Based on the initial condition in the first picture in Figure \subref*{fig:closeoneloopinitialcondinit} we close more loops in the neighborhood of this closed loop in the initial data in Figure \ref{fig:closeloopsleftinitialcond}. Closing iteratively one additional loop results in the initial conditions in Figure \subref*{fig:closeloopsleftinitialcondinit} and the associated steady states are depicted in Figure \subref*{fig:closeloopsleftinitialcondsteadystate}. Note that closing loops close to the source leads to different steady states. In particular, closing loops iteratively in the initial data leads to  steady states which only differ locally. More precisely, the resulting steady states all have the same number of non-zero conductivities. Closing one loop in the initial data results in a steady state which can be  obtained from steady states with the previous initial data by interchanging a non-zero with a negligible conductivity. In particular, the steady states strongly depend on the initial data.

\begin{figure}[htbp]
	\centering
	\subfloat[Initial data] {\includegraphics[width=0.24\textwidth]{microsol_gamma50initialcondnoperturb_poster}\label{fig:initialcondvarynu}}
	\subfloat[$\nu=1$] {\includegraphics[width=0.24\textwidth]{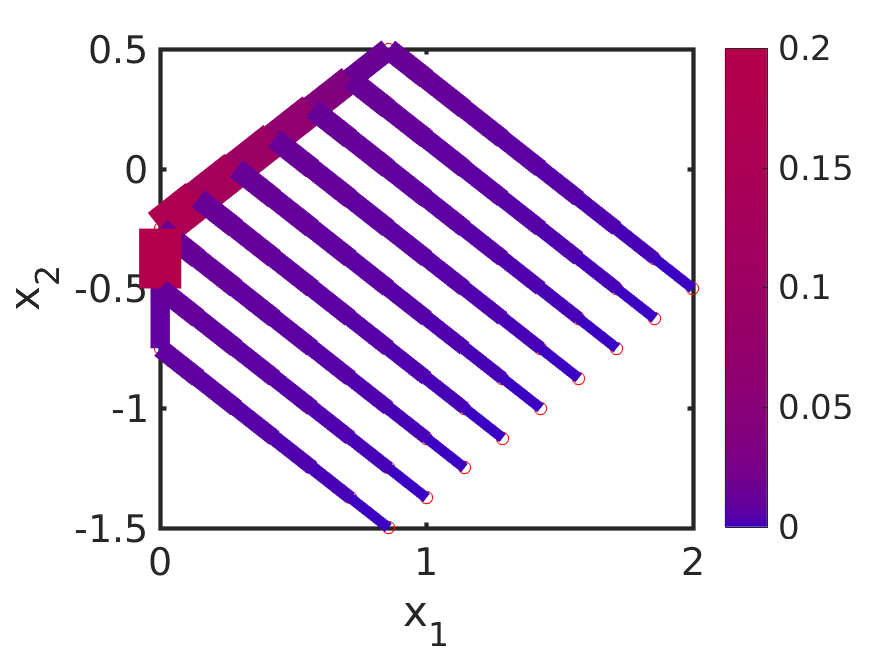}}
	\subfloat[$\nu=100$] {\includegraphics[width=0.24\textwidth]{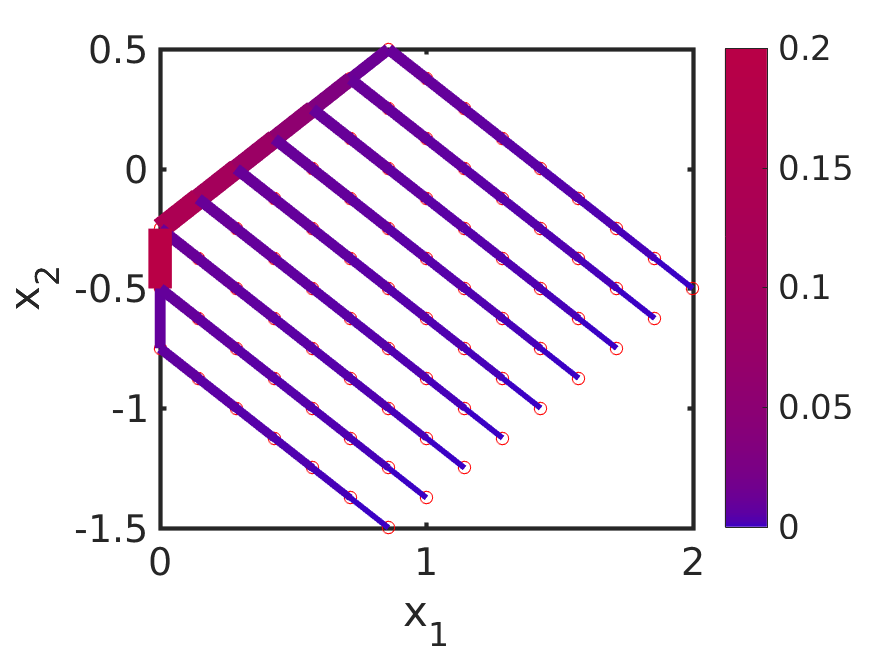}}
	\subfloat[$\nu=10^5$] {\includegraphics[width=0.24\textwidth]{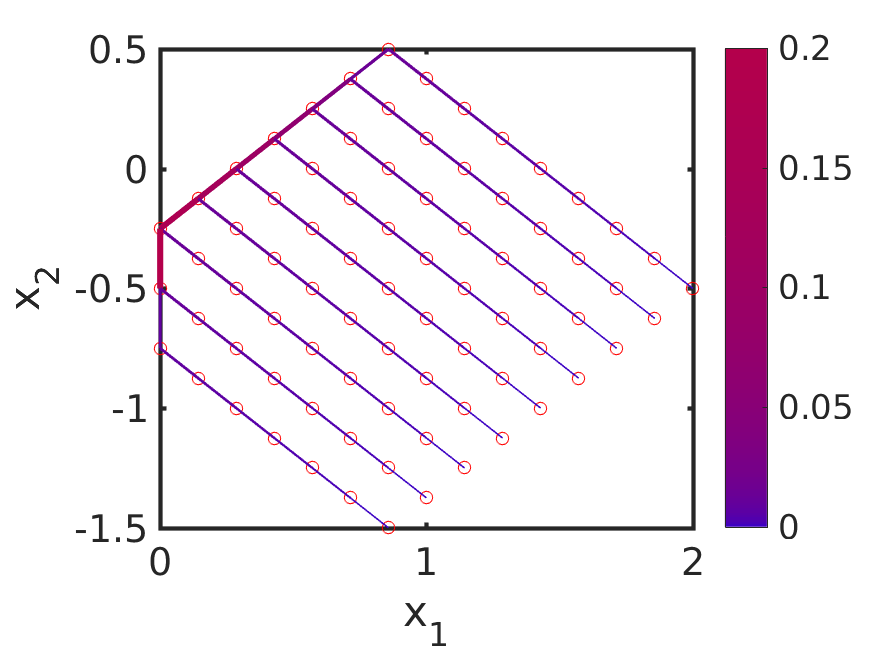}}
	\caption{Steady states  for different values of the parameter $\nu$ in the energy functional \eqref{eq:energydisc}.}\label{fig:varynu}
\end{figure}

In Figure \ref{fig:varynu} the steady states are shown for the same initial data as before (see Figure \subref*{fig:initialcondvarynu})
for different values of the parameter $\nu>0$ in the definition of the energy functional \eqref{eq:energydisc}.
As $\nu$ increases the form of the steady states remain the same, i.e., positive conductivities remain positive for different values of $\nu$. However, the absolute value of the conductivities decreases as $\nu$ increases, see Figure \ref{fig:varynu}. This is consistent with the definition of the energy functional \eqref{eq:energydisc} where the metabolic term is of the form $\frac{\nu}{\gamma} C_{ij}^{\gamma}$ with $\gamma>0$.

\begin{figure}[htbp]
	\centering
	\subfloat[$\delta=100$] {\includegraphics[width=0.24\textwidth]{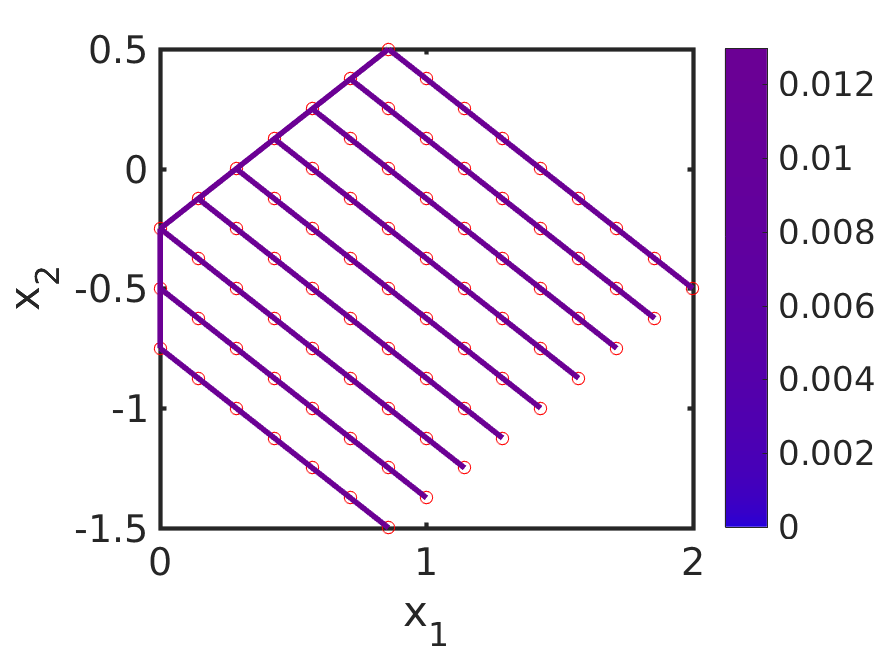}\label{fig:intialcondabsvalueinit1}}
	\subfloat[$\delta=10^3$] {\includegraphics[width=0.24\textwidth]{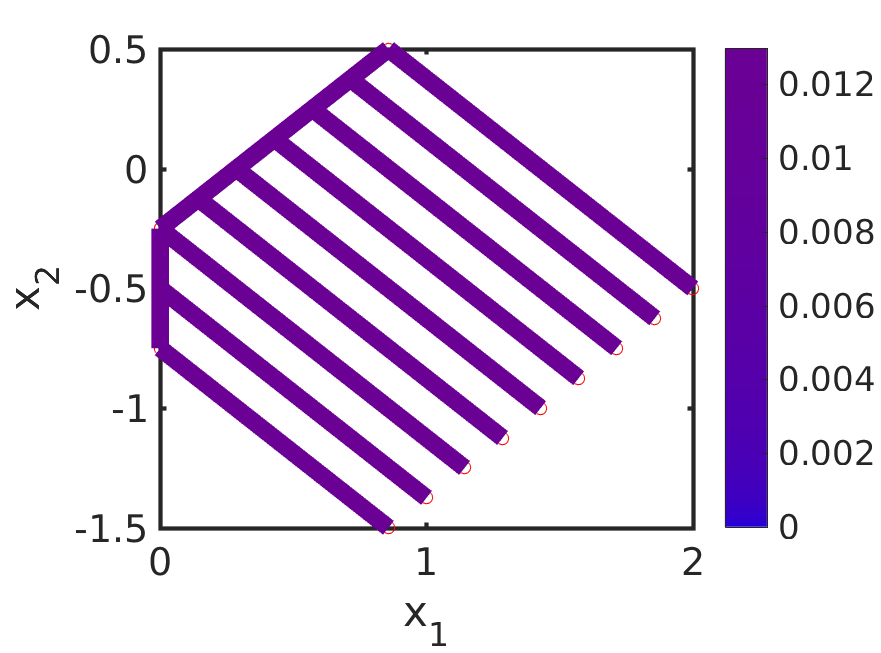}}
	\subfloat[$\delta=10^4$] {\includegraphics[width=0.24\textwidth]{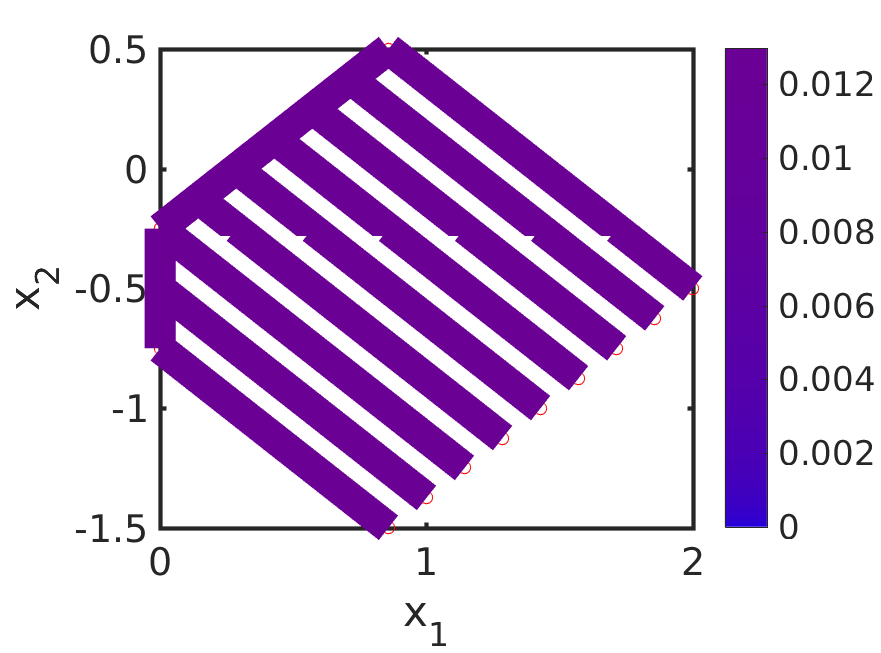}\label{fig:intialcondabsvalueinit3}}
	\subfloat[Steady state] {\includegraphics[width=0.24\textwidth]{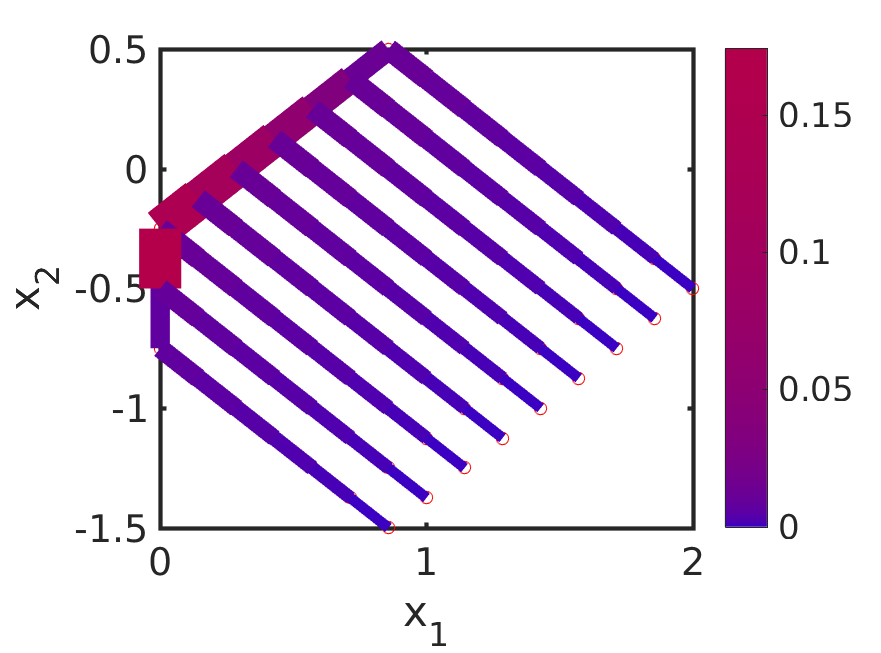}\label{fig:intialcondabsvaluesteadystate}}
	\caption{Initial data for the conductivity vector $\bar{C}=(\bar{C}_{ij})_{(i,j)\in\Eset}$ in the form of a tree where each non-zero conductivity $\bar{C}_{ij}$ is of size $\delta>0$ (left)
	all leading to an identical steady state (right) for the discrete model.}\label{fig:intialcondabsvalue}
\end{figure}

The absolute value of the initial conductivities is varied in Figures \subref*{fig:intialcondabsvalueinit1}--\subref*{fig:intialcondabsvalueinit3} and we show the resulting steady state in Figure \subref*{fig:intialcondabsvaluesteadystate}. More precisely, we consider initial data in the form of a tree as before, when  only those conductivities $\bar{C}=(\bar{C}_{ij})_{(i,j)\in\Eset}$ with positive conductivities $\bar{C}_{ij}$ are considered but vary the absolute value of the initial conductivities. We consider the initial data $\bar{C}_{ij}=\delta$ for every edge $(i,j)\in\Eset$ on the tree for $\delta=5,50,5000,50000$ and $\bar{C}_{ij}=10^{-10}$ otherwise, as shown in Figure \subref*{fig:initialcond} and Figures \subref*{fig:intialcondabsvalueinit1}--\subref*{fig:intialcondabsvalueinit3}, respectively. All these different initial data result in the same steady state  shown in Figure \subref*{fig:intialcondabsvaluesteadystate}. 

\begin{figure}[htbp]
	\centering
	\subfloat[Full graph] {\includegraphics[width=0.24\textwidth]{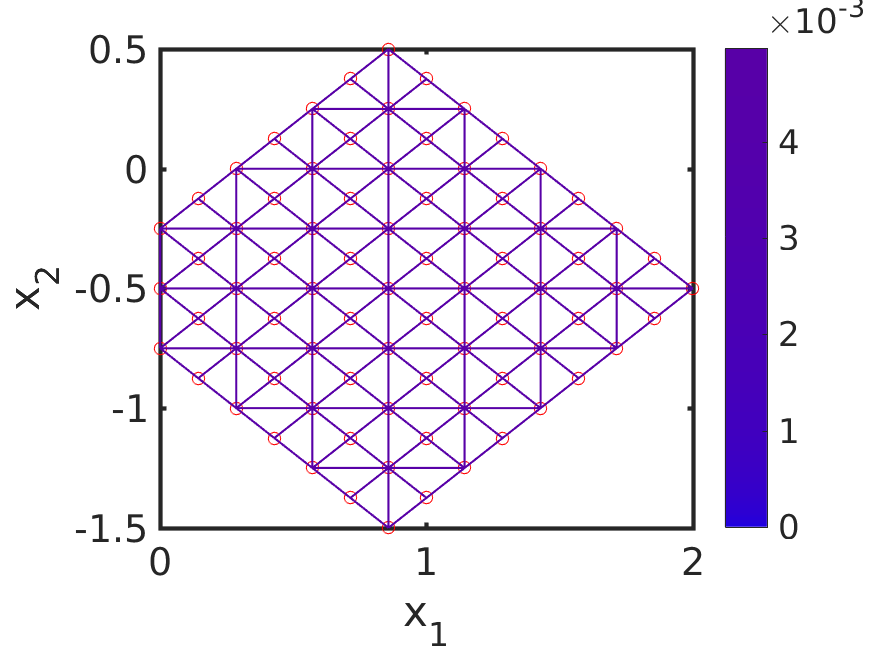}\includegraphics[width=0.24\textwidth]{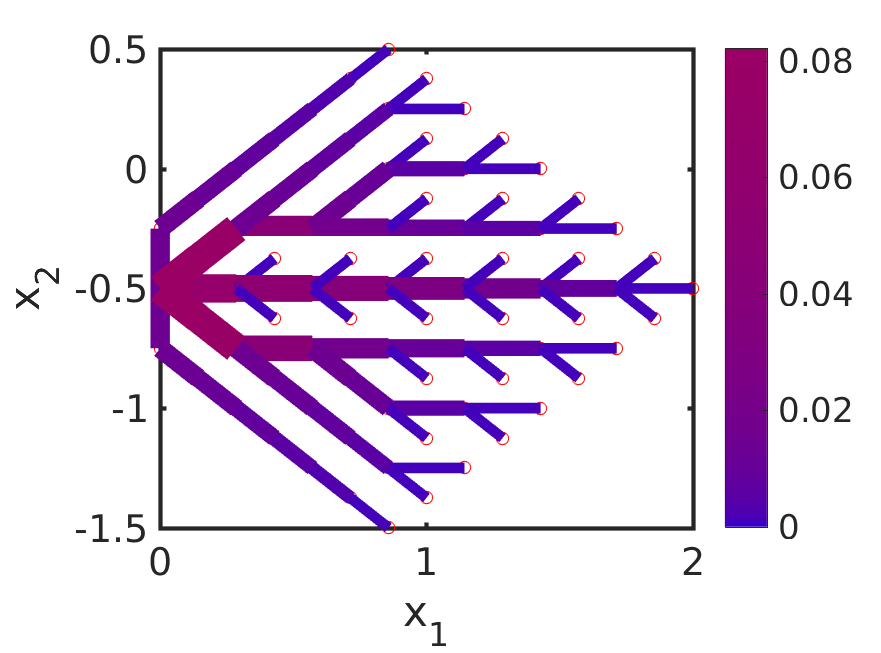}}
	\subfloat[Steady state] {\includegraphics[width=0.24\textwidth]{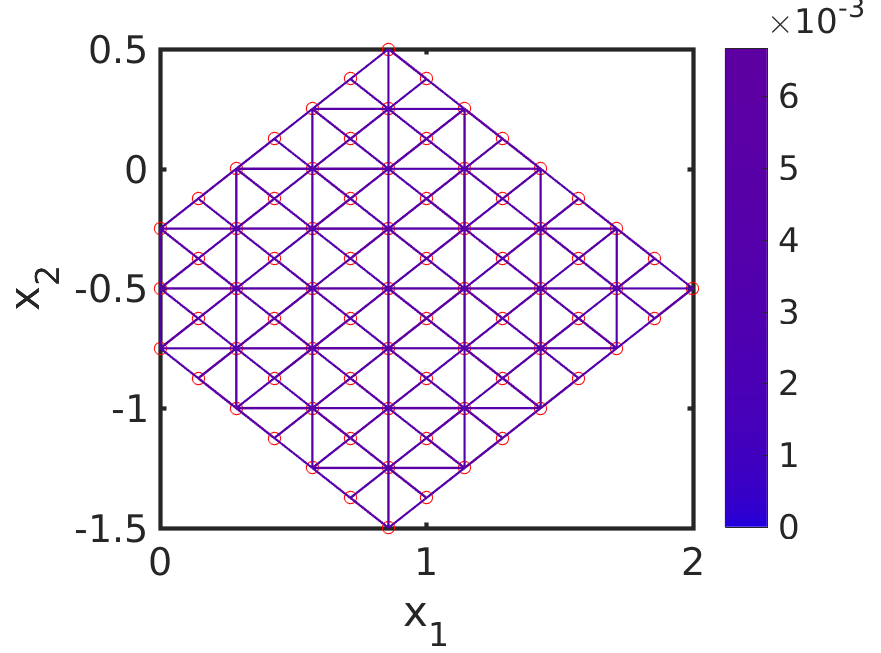}\includegraphics[width=0.24\textwidth]{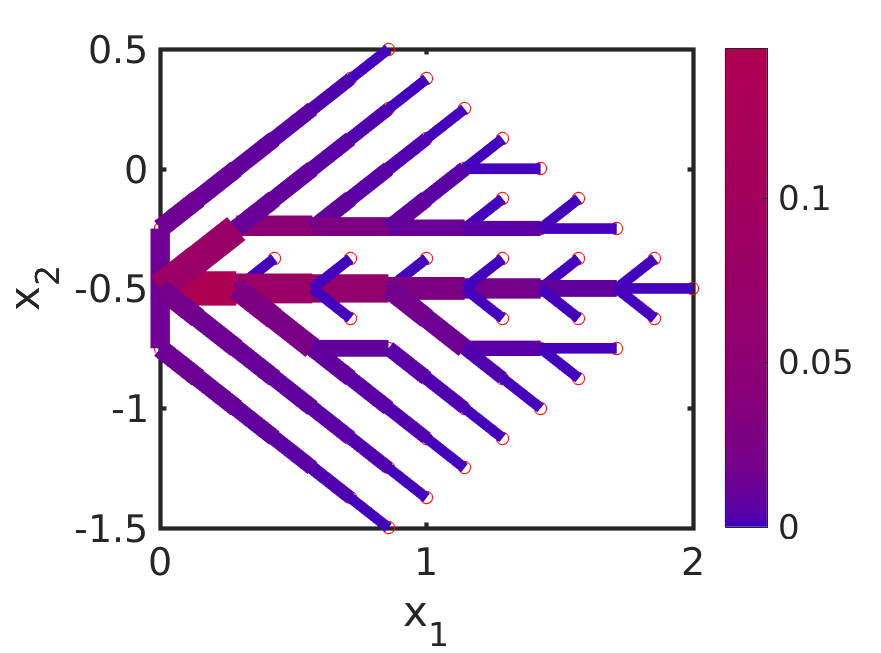}}
	\caption{Steady states  for full graph and perturbed full graph as initial data in the discrete model.}\label{fig:fullgraphinitial}
\end{figure}

In Figure \ref{fig:fullgraphinitial},  full graphs are considered as initial data and we show the associated steady states. We consider  $\bar{C}_{ij}=1$ for all $(i,j)\in\Eset$ and the perturbed full graph with $\bar{C}_{ij}=1+\mathcal{U}(0,1)$ where $\mathcal{U}(0,1)$ denotes a uniformly distributed random variable on $[0,1]$. The associated steady states are more complex transportation networks.

\section*{Acknowledgments}
LMK was supported by the UK Engineering and Physical Sciences Research Council (EPSRC) grant
EP/L016516/1 and the German National Academic Foundation (Studienstiftung des Deutschen Volkes).

\bibliographystyle{plain}
\bibliography{references}
\end{document}